\documentclass{ws-m3as}

\begin{document}




\newcommand{\lip}[1]{W^{#1,\infty}({\mathbb{R}})}
\newcommand{\cb}[1]{{\cal C}^{#1}_b({\mathbb{R}})}
\newcommand{\lipp}[1]{W^{#1,\infty}_L({\mathbb{R}})}
\newcommand{\cbp}[1]{{\cal C}^{#1}_L({\mathbb{R}})}

\newcommand{\cesp}[2]{{\cal C}^{#1}\left(0,T;#2\right)}
\newcommand{\lesp}[2]{L^{#1}\left(0,T;#2\right)}
\newcommand{\wesp}[2]{W^{#1}\left(0,T;#2\right)}
\newcommand{\hesp}[2]{H^{#1}\left(0,T;#2\right)}


\let\Hypothesisfont\itshape
\let\Hypothesisheadfont\bfseries
\newtheorem{hypothesis}{Hypothesis}


\newcommand{\eqref}[1]{{\rm (\ref{#1})}}





\markboth{J. A. Carrillo and S. Labrunie}{Vlasov--Maxwell System for
  Laser-Plasma Interaction}

%
\catchline{}{}{}{}{}
%

\title{GLOBAL SOLUTIONS FOR THE ONE-DIMENSIONAL
VLASOV--MAXWELL SYSTEM FOR LASER-PLASMA INTERACTION}

\author{JOS\'E A. CARRILLO}

\address{ICREA (Instituci\'o Catalana de 
Recerca i Estudis Avan\c{c}ats) \\ and Departament de 
Ma\-te\-m\`a\-ti\-ques \\
Universitat Aut\`onoma de Barcelona \\ 
08193 Bellaterra, Spain.\\
carrillo@mat.uab.es}

\author{SIMON LABRUNIE}

\address{Institut \'Elie Cartan (Math\'ematiques)\\
Universit\'e Henri Poincar\'e Nancy 1 et INRIA (projet CALVI)\\
54506 Vand{\oe}uvre-l\`es-Nancy, France. \\
labrunie@iecn.u-nancy.fr}

\maketitle

\begin{history}
\received{(Day Month Year)}
\revised{(Day Month Year)}
\comby{(xxxxxxxxxx)}
\end{history}

\begin{abstract}
We analyse a reduced 1D Vlasov--Maxwell system introduced recently in 
the physical literature for studying laser-plasma interaction. 
This system can be seen as a standard Vlasov equation in which the 
field is split in two terms: an electrostatic field obtained
from Poisson's equation and a vector potential term satisfying a 
nonlinear wave equation. Both nonlinearities in the Poisson and wave
equations are due to the coupling with the Vlasov equation 
through the charge density. We show global existence of weak solutions 
in the non-relativistic case, and global existence of characteristic 
solutions in the quasi-relativistic case. Moreover, these solutions
are uniquely characterised as fixed points of a certain operator.
We also find a global energy functional for the system allowing us to 
obtain $L^p$-nonlinear stability of some particular equilibria in the 
periodic setting.
\end{abstract}



%
%

\setcounter{equation}{0}


\section{Introduction}
Given a population of electrons, with mass~$m$ and charge~$-e$,
assumed to be relativistic, we denote
$$
\mathbf{v}(\mathbf{p}):= 
\frac{\mathbf{p}}{m\,\sqrt{1+\displaystyle
\frac{|\mathbf{p}|^2}{m^2\,c^2}}} := \frac{\mathbf{p}}{m\,\gamma},
$$
the velocity corresponding to a given momentum~$\mathbf{p}$.
The electrons move under the effect of an electric field~$\mathbf{E}$
and a magnetic field~$\mathbf{B}$.
Then, their distribution function $\mathfrak{f}(t,\mathbf{q},
\mathbf{p})$, where $\mathbf{q}$~denotes the position variable, is 
solution to the Vlasov equation:
\begin{equation}
\frac{\partial\mathfrak{f}}{\partial t} + \mathbf{v}(\mathbf{p})\cdot 
\frac{\partial\mathfrak{f}}{\partial\mathbf{q}} - {e}\, 
\left(\mathbf{E} + \mathbf{v}(\mathbf{p})\times\mathbf{B} \right)\cdot 
\frac{\partial\mathfrak{f}}{\partial\mathbf{p}} = 0.
\label{vla.dim}
\end{equation}
The fields $\mathbf{E}$ and~$\mathbf{B}$ are the sum of three parts:
\begin{enumerate}
\item the self-consistent fields created by the electrons;
\item the electromagnetic field of a laser wave which is sent into
the medium (called the \emph{pump} wave);
\item the electrostatic field $\mathbf{E}_\mathrm{ext}(\mathbf{q})$
generated by a background of ions which are considered immobile during 
the time scale of the wave, and/or by an external, static confinement 
potential.
\end{enumerate}
In all cases, we denote by $n_\mathrm{ext}:=\varepsilon_0\,
\mathrm{div\,}\mathbf{E}_\mathrm{ext}/e$. Without this term,
the population of electrons could not be dynamically stable.
Then, the Maxwell system is written:
\begin{eqnarray}
\frac{\partial\mathbf{E}}{\partial t} &=& c^2\,\mathbf{curl\,B} + 
\frac{e}{\varepsilon_0}\,\mathbf{j},
\label{maxamp.dim}\\
\frac{\partial\mathbf{B}}{\partial t} &=& -\mathbf{curl\,E} ,
\label{fara.dim}\\
\mathrm{div\,}\mathbf{E} &=& \frac{e}{\varepsilon_0}\,
\left(n_\mathrm{ext}-n\right),
\label{gauss.dim}\\
\mathrm{div\,}\mathbf{B} &=& 0,
\label{monop.dim}
\end{eqnarray}
where $c$ and~$\varepsilon_0$ are the speed of light and the dielectric
permittivity of vacuum, and the electron density and flux 
$n$ and~$\mathbf{j}$ are the first two moments of the distribution
function~$\mathfrak{f}$:
\[
\{n,\mathbf{j}\}(t,\mathbf{q}) := 
\int \{1,\mathbf{v}(\mathbf{p})\}\, 
\mathfrak{f}(t,\mathbf{q},\mathbf{p})\, \mathrm{d}\mathbf{p}.
\]
It is well known that Eqs.~(\ref{monop.dim}, \ref{fara.dim})
amount to the existence of vector and scalar potentials such that
\begin{equation}
\mathbf{B} = \mathbf{curl\,A},\quad \mathbf{E} =
- {\partial_t\mathbf{A}} - \mathbf{grad\,}\Phi.
\label{eq:chps:pot}
\end{equation}

The assumptions below the 1D model are the following: all variables
depend on only one space variable, denoted~$x$, and the electrons
are monokinetic in the directions transversal to~$x$. This is 
physically justifed by the fact that all the phenomena, especially 
heating, are much more rapid along the direction of propagation of the 
laser wave than in the transversal directions.
So, the distribution function becomes:
\begin{equation}
\mathfrak{f}(t,\mathbf{q},\mathbf{p}) = f(t,x,p_x)\, \delta\left( 
\mathbf{p}_\bot - \mathbf{p}_0(t,x) \right).
\label{eq:red:dim:f}
\end{equation}

The function $\mathbf{p}_0(t,x)$ can be determined by 
Hamiltonian-mechanical cons\-id\-era\-tions\cite{HGBSC03}. 
The Hamiltonian for one particle is $H:= \gamma\,m\,c^2 - e\,\Phi$, 
and~(\ref{vla.dim}) reads:
\begin{displaymath}
\frac{\partial \mathfrak{f}}{\partial t} + [H,\mathfrak{f}] = 0,
\end{displaymath}
where~$[\cdot,\cdot]$ is the Poisson bracket. Then, by Hamilton's 
equation, the transversal component of the canonical 
conjugate momentum $\mathbf{P_c} := \mathbf{p} -e\,\mathbf{A}$ is 
conserved:
$$\frac{\mathrm{d}\mathbf{P_c}}{\mathrm{d}t} = 
- \frac{\partial H}{\partial\mathbf{q}}
\,\Longrightarrow\, \frac{\mathrm{d}\mathbf{P_{c\bot}}}{\mathrm{d}t} = 
- \frac{\partial H}{\partial\mathbf{q_\bot}} = 0.$$
By a suitable change of referential, we can suppose that 
$\mathbf{P_{c\bot}}=0$; and by imposing the Coulomb gauge
$\mathrm{div\,}\mathbf{A} = 0$, that $A_x=0$. Hence,
$\mathbf{p}_0(t,x) = e\,\mathbf{A}(t,x)$.

For the sake of simplicity, we shall assume in this work that the 
pump wave is linearly polarised in a direction which we call~$y$; 
however, the forthcoming computations can be easily generalised to
an arbitrary polarisation. Under these circumstances, 
Eq.~(\ref{eq:chps:pot}) becomes:
$$B_z = \partial_x A_y,\quad E_y = -\partial_t A_y,\quad 
E_x = -\partial_x\Phi,\quad E_z=B_y=0,$$
which allows to recast the Vlasov equation~(\ref{vla.dim}) and the two 
remaining Maxwell equations~(\ref{maxamp.dim}, \ref{gauss.dim}) 
as the following system:
\begin{eqnarray}
\frac{\partial f}{\partial t} + 
\frac{p_x}{m\gamma}\, \frac{\partial f}{\partial x} - e\,
\left(E_x + \frac{e\,A_y}{m\,\gamma}\,\frac{\partial A_y}{\partial x}
\right)\, \frac{\partial f}{\partial p_x} &=& 0,
\label{vla.dim.2}\\
\frac{\partial^2 A_y}{\partial t^2} - c^2\,
\frac{\partial^2 A_y}{\partial x^2} &=& 
-\frac{e}{m\varepsilon_0}\,n_\gamma\,A_y
\label{ond.dim.2}\\
\frac{\partial E_x}{\partial t} &=& \frac{e}{\varepsilon_0}\,j_x
\label{amp.dim.2}\\
\frac{\partial E_x}{\partial x} &=& \frac{e}{\varepsilon_0}\,
\left(n_\mathrm{ext}-n\right).
\label{gauss.dim.2}
\end{eqnarray}
The Lorentz factor~$\gamma$ and the density, \emph{quasi-density} 
and flux are now given~by:
\begin{eqnarray}
\gamma &=& \sqrt{1+\frac{p_x^2}{m^2\,c^2}+\frac{e^2\,A_y^2}{m^2\,c^2}}
\;;\label{gamma.2}\\
\left\{ n,n_\gamma,j_x \right\} &=& \int \left\{1,\frac1\gamma,
\frac{p_x}{m\gamma}\right\}\, f\,\mathrm{d}p_x.
\label{n.j.2}
\end{eqnarray}

The equations~(\ref{amp.dim.2}) and~(\ref{gauss.dim.2}), which are 
relative to the same variable~$E_x$, are redundant under regularity 
conditions. Eq.~(\ref{amp.dim.2}) is a (simple) evolution equation, 
while~(\ref{gauss.dim.2}) is interpreted as a constraint. As usual, 
the satisfaction of this constraint at $t=0$ implies its satisfaction 
at any time, thanks to~(\ref{amp.dim.2}) and to the continuity
equation $\partial_t n + \partial_x j_x = 0$.

For a general polarisation, (\ref{ond.dim.2}) would be duplicated,
with a similar equation for~$A_z$. Quadratic terms
in $A_y$ and~$A_z$ would be added in~(\ref{gamma.2}) as well as in
the third term in~(\ref{vla.dim.2}). The reader will convince himself 
that the study of this slightly more complicated system is completely
similar to that of~(\ref{vla.dim.2}--\ref{n.j.2}).
{F}rom now on, we shall omit the subscripts in $p_x,\ E_x,\ j_x,\ A_y$.

\medbreak

Let us notice that the model~(\ref{vla.dim.2}--\ref{n.j.2}), as well as
its extended version for a general polarisation, are classes of
\emph{exact} solutions to the relativistic Vlasov--Maxwell model, 
without any approximation. They are, as far as we know, the simplest 
exact solutions beyond 1D~Vlasov--Poisson models.

\subsection{Discussion of the relativistic character}
The model~(\ref{vla.dim.2}--\ref{n.j.2}) features a strongly 
non-linear coupling between the kinetic and electromagnetic variables,
through the Lorentz factor~(\ref{gamma.2}). This phenomenon makes this
system difficult to study on the theoretical level, but also to 
solve numerically: no splitting between the variables is 
possible.\cite{HGBSC03} This is why two reduced models have been 
defined by physicists:
\begin{enumerate}
\item The \emph{non-relativistic model} (hereafter denoted NR) 
approximates the relativistic dynamic by the Newtonian one
by setting $\gamma=1$ everywhere. It is physically justified when 
the temperature is low enough, so that the proportion of 
relativistic electrons is negligible, \emph{and}
the intensity of the pump wave is small.
\item The \emph{quasi-relativistic\footnote{
Also called \emph{semi-relativistic} by some authors.} 
model} (QR) consists in approximating $\gamma$ by 
$\sqrt{1+\frac{p^2}{m^2\,c^2}}$ in the second term 
in~(\ref{vla.dim.2}) and in the definition of~$j$,
and setting $\gamma=1$ in the third term in~(\ref{vla.dim.2}) 
and in the definition of~$n_\gamma$, which
amounts to setting $n_\gamma=n$. It is acceptable when the proportion 
of ultra-relativistic ($v\simeq c$) electrons is negligible and
the pump intensity is moderate.
\end{enumerate}
The original model, with $\gamma$ defined by~(\ref{gamma.2}) will
be referred to as \emph{fully relativistic}~(FR). We remark that 
the NR~model is a class of exact solutions to the non-relativistic 
Vlasov--Maxwell system, \hbox{i.e.}~(\ref{vla.dim}--\ref{monop.dim}) 
with $\mathbf{v}(\mathbf{p})=\mathbf{p}/m$. By contrast, the QR~model 
is only an approximation to the FR~one. What makes its interest for the
applied mathematician --- besides its widespread use for simulation ---
is that it already contains certain features of higher-dimensional 
relativistic Vlasov--Maxwell systems, while being simpler to study.

\subsection{Rescaled equations and Cauchy problem}
The set of equations~(\ref{vla.dim.2}--\ref{n.j.2}) can be simplified
by introducing some rescaled variables. Let $\overline n$ be the unit
of density; we choose the units for the independent variables as:
$$\overline x = \frac{c}{e}\,
\sqrt{\frac{m\,{\varepsilon_0}}{\overline{n}}},\quad 
\overline t = \frac{\overline x}{c},\quad\overline p = m\,c\,;$$
and for the dependent variables as:
$$\overline E = \sqrt{\frac{m\,c^2\,\overline{n}}{{\varepsilon_0}}},
\quad\overline A=\frac{m\,c}e,\quad 
\overline f = \frac{\overline n}{m\,c},
\quad\overline\jmath = \overline n\,c.$$
Keeping the same notations for the rescaled variables, we 
obtain the rescaled system:
\begin{eqnarray}
\frac{\partial f}{\partial t} + 
\frac{p}{\gamma_1}\, \frac{\partial f}{\partial x} -
\left(E + \frac{A}{\gamma_2}\,\frac{\partial A}{\partial x}\right)\, 
\frac{\partial f}{\partial p} &=& 0,
\label{vla.adim}\\
\frac{\partial E}{\partial t} &=& j
\label{amp.adim}\\
\frac{\partial^2 A}{\partial t^2} - \frac{\partial^2 A}{\partial x^2} 
&=& - n_\gamma\,A
\label{ond.adim}
\end{eqnarray}
where the flux~$j$ and the quasi-density~$n_\gamma$ are defined as
$$j(t,x):= \int \frac{p}{\gamma_1}\, f(t,x,p)\,\mathrm{d}p,\qquad
n_\gamma(t,x):= \int \frac{1}{\gamma_2}\, f(t,x,p)\,\mathrm{d}p.$$
Of course, suitable initial conditions are supplied, namely
\begin{equation}
f(0,x,p) = f_0(x,p),\ A(0,x) = A_0(x),\
\partial_t A(0,x) = {\dot A}_0(x).
\label{ci.adim}
\end{equation}
For coherence, the initial electrostatic field must be given by:
\begin{equation}
E_0(x) = E_0(0) - \int_0^x \left(n_0(y) - n_\mathrm{ext}(y)
\right)\, \mathrm{d}y,\quad\mbox{where:}\quad n_0(x) := \int f_0(x,p)
\, \mathrm{d}p.\ 
\label{ci2.adim}
\end{equation}
This guarantees that the Poisson (or Gauss) equation
\begin{equation}
\frac{\partial E}{\partial x} = n_\mathrm{ext} - n,
\quad\mbox{where:}\quad n(t,x):= \int f(t,x,p)\,\mathrm{d}p,
\label{poi.adim}
\end{equation}
will hold at any time provided the continuity equation is satisfied. 

As far as the relativistic character of the particles is concerned,
the three versions of the model are respectively:
\begin{itemize}
\item NR: $\gamma_1=\gamma_2=1$.
\item QR: $\gamma_1=\sqrt{1+p^2},\ \gamma_2=1$.
\item FR: $\gamma_1=\gamma_2=\sqrt{1+p^2+A^2}$.
\end{itemize}

Let us note that we shall investigate the existence of two classes 
of solutions for the system~(\ref{vla.adim}--\ref{ond.adim}): 
periodic solutions, corresponding to initial data 
that are periodic in space, with a given period $L$, and 
``open-space'' solutions, \hbox{i.e.}~solutions of finite mass and 
energy. In both cases, we always assume that $n_\mathrm{ext}$ is at 
least bounded.
Moreover, in the periodic setting, we assume that it is periodic and
$$\int_0^L \left(n_\mathrm{ext}(x) - n_0(x)\right)\,\mathrm{d}x =0,$$
so that $E_0$ is indeed periodic. 

\medbreak

To the best of our knowledge, this is the first mathematical work on
this particular Vlasov--Maxwell 
system~(\ref{vla.adim}--\ref{ond.adim}): most of the previous 
mathematical works on reduced Vlasov--Maxwell models deal with 
systems living in two or one-and-a-half 
dimensions.\cite{GlSt,GlSc,Glassey96,Bouchut} In our case, the
lower dimension is, to some extent, compensated by a stronger 
nonlinearity. This system shares some
common features with the Nordstr\"{o}m--Vlasov system recently studied
in \cite{CalRe,CalRe2}. Both systems present a Vlasov equation coupled
to a wave equation whose right-hand side depends on the charge
density. The main differences between them are the gravitational 
character of the Nordstr\"{o}m--Vlasov system, our coupling with the
Poisson equation and our more complicated right-hand side in the 
wave equation.

In this article, we shall only investigate the existence 
and uniqueness of solutions to the NR and QR models.
The differences between them can be bridged easily by
using the shorthand notation
$$
\widehat p:= p \mbox{ (NR case)},\quad
\widehat p:= \frac{p}{\sqrt{1+p^2}} \mbox{ (QR case)}.
$$
On the other hand, the method presented in this article cannot apply
directly to the FR~model, because of its much stronger and more
non-linear coupling between the kinetic and electromagnetic variables.

\medbreak

The next Section is devoted to reviewing some basic estimates for the
Vlasov equation, and the different notions of solutions we will 
deal with. In Section~3 we use the procedure of~\cite{CoKl80} to
prove the global existence of weak
solutions in the NR and QR~cases, which are characterised as 
unique fixed points of a certain operator. 

This result is improved in two ways in Sections 4 and~5
respectively. We first prove the local-in-time existence of 
characteristic solutions in the NR case. The main difficulty in 
proving global existence of characteristic solutions consists in
controlling the second space derivative of the vector potential $A$ 
or, equivalently, the first space derivative of the density~$n$. These
difficulties are analogous to those found in the two and one-half 
dimensional Vlasov--Maxwell systems studied 
in \cite{GlSt,GlSc,Glassey96,Bouchut}.
In the QR~case, we obtain the existence result thanks to a good 
integral representation of the second derivative of $A$ reminiscent of 
similar ideas in \cite{GlSc}, used recently in \cite{CalRe2}.

Finally, the last section is devoted to obtaining a global energy
functional for the three cases. It is conserved in time by 
characteristic solutions in the NR and QR cases.
This fact, together with now standard relative entropy
arguments, leads to the $L^p$-nonlinear stability of a family of
steady states in the periodic setting. 


\section{Solutions to the forced Vlasov equation}
In this Section, we assume that the fields $(E,A)$ are given
and we introduce several notions of solution to~(\ref{vla.adim}) and 
summarise their associated regularity properties.

Whatever their regularity, we are mainly interested in \emph{global} 
solutions, \hbox{i.e.}~which exist for any time. However, our 
estimates will generally not be uniform in time and thus, we fix an 
\emph{arbitrary} 
target time~$T$ and we look for solutions defined on~$[0,T]$, for all 
$T>0$. We use the following notations for functional spaces:
\begin{itemize}
\item $\cesp{m}{{\cal X}}$: the space of $m$~times continuously 
differentiable functions from $[0,T]$ to the Banach space~${\cal X}$;
\item $\lip{k}$: the space of functions from ${\mathbb{R}}$ 
to~${\mathbb{R}}$ with all derivatives (in the sense of distributions) 
bounded, up to order~$k$;
\item $\cb{k}={\cal C}^k({\mathbb{R}})\cap\lip{k}$: the space of 
$k$~times 
continuously differentiable functions with all derivatives bounded;
\item $\lipp{k},\ \cbp{k}$: the subspaces of functions which have the
space period~$L$.
\end{itemize}
We will denote by $\left\|\cdot\right\|_t$ the norm
in $L^\infty((0,t)\times{\mathbb{R}})$, for~$t\le T$.

\subsection{The characteristic system}
{F}rom now on, we denote
\begin{equation}
F(t,x):= E(t,x) + A(t,x)\,\partial_x A(t,x)
\label{df:force}
\end{equation}
the force generated by the fields $(E,A)$.
The characteristic system associated to the transport 
equation~(\ref{vla.adim}) reads:
\begin{equation}
\left\lbrace
\begin{array}{ll}
{\displaystyle \frac{\mathrm{d}X}{\mathrm{d}s}} = \widehat P(s),\quad&
{\displaystyle \frac{\mathrm{d}P}{\mathrm{d}s}} = -F(s,X(s)),\hfill \\
X(t)=x,\quad& P(t)=p.\hfill
\end{array}
\right.
\label{carac}
\end{equation}
Global existence and uniqueness of solution to the above 
system is ensured by assuming that the force field is continuous 
in time and globally Lipschitz in space; in turn, a sufficient 
condition for this is:
\begin{equation}
E\in\cesp{0}{\lip{1}},\quad A\in\cesp{0}{\lip{2}}.
\label{lip-reg}
\end{equation}
Under this assumption, the unique solution to the characteristic 
system~(\ref{carac}) denoted by $(X(s;t,x,p),P(s;t,x,p))$ 
becomes (at least) a continuous function in all its variables. 
We shall also consider a stronger regularity condition
\begin{equation}
E\in\cesp{0}{\cb{1}},\quad A\in\cesp{0}{\cb{2}}.
\label{cb-reg}
\end{equation}

{F}rom the uniqueness of the solution, we deduce a periodicity result:

\begin{lemma}
If the force field is periodic in space, \hbox{i.e.}: $\forall (t,x),\ 
F(t,x+L) = F(t,x)$, then the following identity holds for all 
$s,\ t,\ x,\ p$:
\begin{equation}
X(s;t,x+L,p) = X(s;t,x,p) + L,\quad P(s;t,x+L,p) = P(s;t,x,p).
\label{eq:perio:X.P}
\end{equation}
\label{lem:perio:X.P}
\end{lemma}

The divergence of characteristics generated by different force fields 
is measured in a classical way (see for instance \cite[Lemma~1]{CoKl80}
or \cite{Bouchut} and references therein).

\begin{lemma} 
{L}et $(X^1,P^1)$ and~$(X^2,P^2)$ be the characteristics associated
to the respective forces $F_1,\ F_2 \in \cesp{0}{\lip{1}}$. 
Then, the following inequalities hold for all $(t,x,p)$:
\begin{eqnarray}
\left| X^1(0;t,x,p) - X^2(0;t,x,p)\right| &\le& t\,
\int_0^t \left\|F_1 - F_2\right\|_s\, \mathrm{d}s,
\label{42.2}\\
\left| P^1(0;t,x,p) - P^2(0;t,x,p)\right| &\le&
\int_0^t \left\|F_1 - F_2\right\|_s\, \mathrm{d}s.
\label{42.4}
\end{eqnarray}
\end{lemma}
\begin{proof}
Eq.~(\ref{42.4}) clearly stems from the integration of the equation 
for~$P$ in~(\ref{carac}). In the NR~case, Eq.~(\ref{42.2}) immediately 
follows. In the QR~case, we notice that
$$\left| \frac{\mathrm{d}}{\mathrm{d}p} \left[ \frac{p}{\sqrt{1+p^2}} 
\right] \right|
= \frac{1}{\left(1+p^2\right)^{3/2}} \le 1,$$
so $\left| \widehat{P^1} - \widehat{P^2}\right|\le
\left|{P^1} -{P^2}\right|$, hence~(\ref{42.2}) by integration.
\end{proof}

In a very similar fashion, the derivatives of the characteristics can 
be estimated in terms of derivatives of the force field.

\begin{lemma}
The solution $(X,P)$ to~\eqref{carac}, with
$F\in\cesp{0}{\lip{1}}$, is Lipschitz \hbox{w.r.t.} the
variables~$(x,p)$, and thus a.e. differentiable. The derivatives 
satisfy the following bound for a.e. $(x,p)$ and all $0\le\tau\le t$:
\begin{equation}
\max \left\{
\left| \frac{\partial X}{\partial x}(\tau) \right|,\ 
\left| \frac{\partial X}{\partial p}(\tau) \right|,\
\left| \frac{\partial P}{\partial x}(\tau) \right|,\ 
\left| \frac{\partial P}{\partial p}(\tau) \right| \right\} 
\le \mathrm{e}^{(t-\tau)\, \left(1+\left\|\partial_x F\right\|_t
\right)}.
\label{eq:char:der}
\end{equation}
\label{lem:char:der}
\end{lemma}
\begin{proof}
Consider two final conditions~$(x_1,p_1)$ and~$(x_2,p_2)$.
We use the shorthand:
$$\left(X_i(s),P_i(s)\right):= 
\left(X(s;t,x_i,p_i),P(s;t,x_i,p_i)\right),\quad i=1,\ 2 \,;$$
and we denote by $\Lambda F(s)$ the Lipschitz constant of $F(s,\cdot)$,
\hbox{i.e.}~$\sup_x \left| \partial_x F(s,x) \right|$.

By integrating the characteristic system~\eqref{carac}, we get:
\begin{eqnarray*}
&& \left| X_1(\tau) - X_1(\tau) \right| +
\left| P_1(\tau) - P_2(\tau) \right| \\
&&\mbox{}\;\le
 |x_1 - x_2| + |p_1 - p_2| + \int_\tau^t \left\{ 
\left| \widehat P_1(s) - \widehat P_2(s) \right|
+ \left| F(s,X_1(s)) - F(s,X_2(s)) \right| \right\}\,\mathrm{d}s \\
&&\mbox{}\;\le 
 |x_1 - x_2| + |p_1 - p_2| + \int_\tau^t \left\{ 
\left| P_1(s) - P_2(s) \right| + \Lambda F(s)\,
\left| X_1(s) - X_2(s) \right| \right\}\,\mathrm{d}s \\
&&\mbox{}\;\le 
 |x_1 - x_2| + |p_1 - p_2| + \left(1+\left\|\partial_x F\right\|_t
\right)\! \int_\tau^t\! \left\{  \left| X_1(s) - X_2(s) \right| + 
\left| P_1(s) - P_2(s) \right| \right\}\,\mathrm{d}s.
\end{eqnarray*}
Hence, by Gronwall's lemma:
$$ \left| X_1(\tau) - X_1(\tau) \right| + 
\left| P_1(\tau) - P_2(\tau) \right| \le \left(
 |x_1 - x_2| + |p_1 - p_2| \right)\, 
\mathrm{e}^{(t-\tau)\, \left(1+\left\|\partial_x F\right\|_t\right)},$$
which proves the Lipschitz character of the functions $X,\ P$ and the
quantitative estimate~\eqref{eq:char:der}.
\end{proof}

\subsection{Characteristic, mild and weak solutions}
In this work, we shall always assume the following two 
hypotheses about the initial distribution function~$f_0$.

\begin{hypothesis}
In the open-space case, 
$f_0\in L^1({\mathbb{R}}^2)$; in the periodic case,
\begin{eqnarray*}
f_0\in L^1_L({\mathbb{R}}^2)&:=& 
\biggl\{ f\in L^1_\mathrm{loc}({\mathbb{R}}^2) :
\ f(x+L,p) = f(x,p) \, \mbox{a.e. in } (x,p)\in{\mathbb{R}}^2\\
&&\hphantom{\biggl\lbrace} \mbox{ and }
\left. f \right|_{(0,L)\times{\mathbb{R}}} \in L^1((0,L)\times
{\mathbb{R}})\biggr\}.
\end{eqnarray*}
In both cases, $(\partial_x f_0,\partial_p f_0)\in L^1_\mathrm{loc}
({\mathbb{R}}^2)$.
\label{hypo:f0c1}
\end{hypothesis}
\smallbreak
\begin{hypothesis}
There exists a continuous, positive, even function $g(p)$, 
which moreover is decreasing in~$|p|$ and satisfies
$$\int_{\mathbb{R}} |p|\, g(p)\,\mathrm{d}p < \infty,$$
such that
$$f_0(x,p)\le g(p),\quad \left|\partial_x f_0(x,p)\right| \le g(p),
\quad \left|\partial_p f_0(x,p) \right|\le g(p).$$
\label{hypo:f0}
\end{hypothesis}
Let us notice that these hypotheses imply, in the open-space case, that
$f_0\in W^{1,1}_\mathrm{loc}({\mathbb{R}}^2)\cap 
W^{1,\infty}({\mathbb{R}}^2)$; the initial density 
$n_0\in L^1({\mathbb{R}})\cap\lip{1}$; and the initial electrostatic 
field, given by \eqref{gauss.dim.2}, $E_0\in\lip{1}$ if 
$n_\mathrm{ext}\in L^\infty({\mathbb{R}})$, $E_0\in\cb{1}$ if 
$n_\mathrm{ext}\in\cb{0}$. 
In the periodic case, one has similar properties 
in terms of periodic spaces; indeed, electrical neutrality ensures 
that $n-n_\mathrm{ext}$ admits a periodic primitive.

For further reference, let us note that the function $g_r(p)$ defined 
as
$$g_r(p) = g(0) \mbox{ for } |p|\le r,\quad
g_r(p) = g(|p|-r) \mbox{ for } |p|\ge r,$$
satisfies
\begin{equation}
\int_{\mathbb{R}} g_r(p)\, \mathrm{d}p = 2r\,g(0) + \int_{\mathbb{R}} 
g(p)\,\mathrm{d}p.
\label{eq:int:grp}
\end{equation}

The first notion of solutions which we consider is the one given
by the usual characteristic method.

\begin{definition}
Given a force field $F\in\cesp{0}{\lip{1}}$ and $f_0$ satisfying
Hypotheses~\mbox{\rm\ref{hypo:f0c1}--\ref{hypo:f0}},
we define the \emph{characteristic} solution of the Vlasov 
equation~\eqref{vla.adim} as:
\begin{equation}
f(t,x,p):= f_0(X(0;t,x,p),P(0;t,x,p)),
\label{eq:df:char}
\end{equation}
where $(X(s;t,x,p),P(s;t,x,p))$ is the unique solution 
to~\eqref{carac}. Moreover, assuming $F\in\cesp{0}{\cb{1}}$
and $f_0\in{\cal C}^1({\mathbb{R}}^2)$, then we shall often refer to 
the characteristic solution as a \emph{classical} solution.
\end{definition}

The previous definition gets clarified by the following result.

\begin{lemma}
Under the assumption~\eqref{lip-reg}, the characteristic solution
$f$ belongs to $\wesp{1,\infty}{W^{1,\infty}({\mathbb{R}}^2)}$.
If, moreover, $f_0\in{\cal C}^1({\mathbb{R}}^2)$ and~\eqref{cb-reg} 
holds, then $f\in\cesp{1}{{\cal C}^{1}_b({\mathbb{R}}^2)}$ 
and~\eqref{vla.adim} is satisfied in the classical sense.
\label{lem:cont:f}
\end{lemma}
\begin{proof}
According to classical dynamical system theory,
$(X(s;t,x,p),P(s;t,x,p))$ is ${\cal C}^1$ \hbox{w.r.t.}~$s$, and also 
\hbox{w.r.t.}~$t$ given the symmetry of these variables. The Lipschitz
character in~$(x,p)$ was obtained in Lemma~\ref{lem:char:der}. 
The first part of the conclusion then follows from 
Hypothesis~\ref{hypo:f0}, and the fact that the composition of
Lipschitz functions is Lipschitz.

Moreover, under~(\ref{cb-reg}), the force is~${\cal C}^1$; as 
$p\mapsto\widehat p$ is~${\cal C}^\infty$, we deduce that the solution
to~\eqref{carac} is~${\cal C}^1$. The last statement then follows from 
the chain rule.
\end{proof}

Mild solutions are introduced for relaxing the assumption of 
differentiability of the characteristics but keeping the fact that 
they define a family of Lipschitz homeomorphisms in phase space 
(see \cite{Bostan} and references therein). 
In fact, this definition can be rephrased 
using the concept of push-forward of a density through a map, which is 
quite well-known in mass transport theory.\cite{Villani} The 
pushed-forward measure 
$T_\#\rho$ of a given measure $\rho$ in ${\mathbb{R}}^n$ assigns mass
$$T_\#\rho [K] := \rho[T^{-1}(K)]$$
to each Borel set $K \subset {\mathbb{R}}^n$. By this property, it 
satisfies that
\[ 
\int_{{\mathbb{R}}^n} \psi\, \mathrm{d}(T_\#\rho) = 
\int_{{\mathbb{R}}^n} (\psi\circ T)\, \mathrm{d}\rho 
\]
for all test functions $\psi\in {\cal C}^0_c({\mathbb{R}}^n)$.

\begin{definition}
Given a force field $F\in\cesp{0}{\lip{1}}$ and 
$f_0\in L^1_\mathrm{loc}({\mathbb{R}}^2)$, we say that a weakly 
continuous function $f(t,x,p)\in {\cal C}_w\left([0,T];
L^1_\mathrm{loc}({\mathbb{R}}^2) \right)$ is a \emph{mild} solution 
of the Vlasov equation~\eqref{vla.adim} if it satisfies
\begin{equation}
\int_{{\mathbb{R}}^2} \psi(x,p)\,f(t,x,p)\, \mathrm{d}x\,\mathrm{d}p = 
\int_{{\mathbb{R}}^2} \psi(X(t;0,x,p),P(t;0,x,p))\,f_0(x,p)\, 
\mathrm{d}x\,\mathrm{d}p
\label{eq:df:mild}
\end{equation}
for all test functions $\psi\in {\cal C}^0_c({\mathbb{R}}^2)$ and all 
$t\ge0$, \hbox{i.e.}, 
$$f(t,x,p) = (X(t;0,x,p),P(t;0,x,p))_\# f_0$$
for all $t\ge0$.
\end{definition}

Taking into account the change of variables formula for Lipschitz 
functions,\cite{Evans} we deduce that a characteristic solution is a 
mild solution of the Vlasov equation~\eqref{vla.adim}.

We can relax even more the assumptions on the force field and talk
about distributional solutions for the Vlasov 
equation~\eqref{vla.adim}.

\begin{definition}
Given a force field $F\in L^\infty((0,T)\times{\mathbb{R}})$ and 
$f_0\in L^1_\mathrm{loc}({\mathbb{R}}^2)$, we say
that $f(t,x,p)\in L^1_\mathrm{loc}((0,T)\times{\mathbb{R}}^2)$ is a 
\emph{distributional} 
solution of the Vlasov equation~\eqref{vla.adim} if it satisfies
\begin{equation}
-\int_0^T\!\!\! \int_{{\mathbb{R}}^2} \!\left( 
\frac{\partial \psi}{\partial t} + \widehat p\, 
\frac{\partial \psi}{\partial x} - F \, 
\frac{\partial \psi}{\partial p} \right) \! f \, 
\mathrm{d}x\, \mathrm{d}p\,\mathrm{d}t\! 
= \! \int_{{\mathbb{R}}^2} \psi(0,x,p)\,f_0\,\mathrm{d}x \,\mathrm{d}p
\label{eq:df:weak}
\end{equation}
for all test functions $\psi\in {\cal C}_c^\infty([0,T)\times
{\mathbb{R}}^2)$.
\end{definition}

It is easy to check that any mild solution is a distributional 
solution. Moreover, in the particular case of characteristic solutions,
in which $f_0\in W^{1,\infty}({\mathbb{R}}^2)$, $f$ belongs to 
$\wesp{1,\infty}{W^{1,\infty}({\mathbb{R}}^2)}$ by 
Lemma~\ref{lem:cont:f}. 
Therefore, one can check from the weak formulation \eqref{eq:df:weak} 
that $f$ satisfies the Vlasov equation~\eqref{vla.adim} as an
equality almost everywhere of locally bounded functions on 
$(0,\infty)\times{\mathbb{R}}^2$.

\subsection{\emph{A priori} estimates}
Since bounds for functions in $\cb{1}$ or~$\lip{1}$ can often be 
obtained in the same manner, we shall treat in the following 
the two types of characteristic solutions (classical and mild) 
together.

Thanks again to Hypothesis~\ref{hypo:f0}, one can define
the density and flux:
\begin{equation}
\{n,j\}(t,x) := \int_{\mathbb{R}} \{1,\widehat p\}\, f(t,x,p)\, 
\mathrm{d}p.
\label{df:dens:cour}
\end{equation}
Indeed, we have the following more general lemma for moments of 
characteristic solutions of (\ref{vla.adim}). Let us denote by
$m_k (t,x)$ the moment of order~$k\in{\mathbb{N}}$ of a solution 
$f(t,x,p)$ with a given force field $F$.

\begin{lemma}
Let $f$ be a characteristic solution of the Vlasov
equation~\eqref{vla.adim} with force field $F \in \cesp{0}{\lip{1}}$.
Assume that the function $g(p)$ from 
Hypothesis~\mbox{\rm\ref{hypo:f0}} 
has the moment of order k bounded, \hbox{i.e.},
$$M_k := \int_{\mathbb{R}} |p|^k\, g(p)\,\mathrm{d}p < \infty,$$
then $m_k (t,x)$ is well-defined for all $(t,x)$ and
\begin{equation}
\left|m_k (t,x) \right| \le M_k + R_k\left(M_0+M_k, 
t \left\|F\right\|_t \right)
\label{eq:mk} 
\end{equation}
where $R_k(a,b)$ will be defined below. In particular, for $k=0$ we 
obtain
\begin{equation}
\|n\|_t \le M_0 + 2\,g(0)\,t\, \|F\|_t .
\label{44.3}
\end{equation}
\label{lem:moments}
\end{lemma}
\begin{proof}
{F}rom~\eqref{eq:df:char} and~Hypothesis~\ref{hypo:f0}, we deduce:
$$|p|^k f(t,x,p) \le |p|^k g\left(\left|P(0;t,x,p)\right|\right).$$
Now, Eq.~(\ref{42.4}) with one of the force fields replaced by~$0$ 
yields $\left|P(0;t,x,p) - p\right| \le t\,\left\|F \right\|_t$;
as $g$~is decreasing, this gives:
$g\left(\left|P(0;t,x,p)\right|\right) \le g_{t\,\left\|F \right\|_t}
(p)$,
and thus
$$
|p|^k f(t,x,p) \le |p|^k g_{t\,\left\|F\right\|_t}(p)
$$
which is clearly integrable in $p$. This proves that $m_k(t,x)$ is
well-defined and
\begin{eqnarray*}
\left|m_k (t,x)\right| \le \int |p|^k g_{t\,\left\|F\right\|_t}(p)\, 
\mathrm{d}p.
\end{eqnarray*}
Moreover, we notice that for $k\geq 1$ we have
\begin{eqnarray}
\int_{\mathbb{R}} |p|^k\,g_r(p)\,\mathrm{d}p & = & \frac{2\,g(0)}{k+1}
\,r^{k+1} + 
\sum_{i=1}^k C_k^i\,r^{k-i}\,
\int_{\mathbb{R}} |p|^i\,g(p)\,\mathrm{d}p \nonumber \\
&\leq & \frac{2\,g(0)}{k+1}\,r^{k+1} + \left(\sum_{i=1}^k
C_k^i\,r^{k-i}\right)\, \int_{\mathbb{R}} (1+|p|^k)\,g(p)\,\mathrm{d}p 
\nonumber \\
&\leq & \frac{2\,g(0)}{k+1}\,r^{k+1} + \left(\sum_{i=1}^k
C_k^i\,r^{k-i}\right)\, (M_0+M_k) \nonumber \\
& := & R_k(M_0+M_k,r).
\end{eqnarray}
Finally, combining previous inequalities, we obtain \eqref{eq:mk}. 
The estimate on the density \eqref{44.3} follows directly from
previous arguments and \eqref{eq:int:grp} which defines the 
function~$R_0$.
\end{proof}

We can also estimate the divergence of moments corresponding to two
different solutions of the Vlasov equation~(\ref{vla.adim}).

\begin{lemma}
Let $f_1,\ f_2$ be characteristic solutions of the Vlasov 
equation~\eqref{vla.adim} with forces $F_1,\ F_2 \in \cesp{0}{\lip{1}}$
respectively. Assume that the function $g(p)$ from 
Hypothesis~\mbox{\rm\ref{hypo:f0}} has the moment of order~$k$ 
bounded, then 
\begin{equation}
\left|m_{1,k} - m_{2,k}\right|(t,x) \le 
R_k\left(M_0+M_k, \max_i \left\|F_i\right\|_t \right)\,
\int_0^t \left\|F_1 - F_2\right\|_s\, \mathrm{d}s .
\label{eq:mk1-mk2}
\end{equation}
\label{lem:difmoments}
\end{lemma}
\begin{proof}
To estimate $m_{1,k}-m_{2,k}$, we consider the characteristics
$(X^1,P^1)$ and $(X^2,P^2)$ associated to $F_1$ and
$F_2$ respectively. Using the shorthand notation for the 
characteristics
$\left(X^i_0,P^i_0\right) = \left(X^i(0;t,x,p),P^i(0;t,x,p)\right)$,
we write:
\begin{eqnarray*}
\left(m_{1,k} - m_{2,k}\right)(t,x) &=& \int_{\mathbb{R}} |p|^k \left[ 
f_0\left(X^1_0,P^1_0\right) - f_0\left(X^2_0,P^2_0\right) 
\right]\, \mathrm{d}p
\end{eqnarray*}
and thus,
\begin{eqnarray*}
\left|m_{1,k} - m_{2,k}\right|(t,x) \le \int_{\mathbb{R}} |p|^k & &
\left\{ \left|\frac{\partial f_0}{\partial x_0}\left(\tilde X,P^1_0
\right)\right|\, \left| X^1_0 - X^2_0 \right|  \right. \\ & &
+\left. \left|\frac{\partial f_0}{\partial p_0}\left(X^2_0,\tilde P
\right)\right| \, \left| P^1_0 - P^2_0 \right| \right\}\, \mathrm{d}p,
\end{eqnarray*}
where we have made use twice of the one-dimensional Taylor--Lagrange
formula; $\tilde X$, respectively $\tilde P$, lie between
$X^1_0$ and~$X^2_0$, resp.~$P^1_0$ and~$P^2_0$. Then we 
invoke~(\ref{42.2}--\ref{42.4}) to bound
\begin{equation}
\left|m_{1,k} - m_{2,k}\right|(t,x) \le 
\int_{\mathbb{R}} |p|^k \left\{ g\left(\left|P^1_0\right|\right) +
t\, g\left(|\tilde P|\right) \right\}\, \mathrm{d}p \times
\int_0^t \left\|F_1 - F_2\right\|_s\, \mathrm{d}s.
\label{eq:m1k-m2k:1}
\end{equation}
Applying again Eq.~(\ref{42.4}) with one of the force fields replaced 
by~$0$ yields $\left|P^i_0 - p\right| \le t\,\left\|F_i\right\|_t$;
as $g$~is decreasing, this gives:
$g\left(\left|P^i_0\right|\right) \le g_{t\,\left\|F_i\right\|_t}(p)$,
and then:
$$g\left(|\tilde P|\right) \le \max\left\{
g\left(\left|P^1_0\right|\right) , g\left(\left|P^2_0\right|\right)
\right\} \le \max_i g_{t\,\left\|F_i\right\|_t}(p).$$
Moreover, $r\mapsto g_r(p)$ is an increasing function of~$r$, for 
all~$p$; this implies:
\begin{eqnarray*}
\int_{\mathbb{R}} |p|^k \left\{ g\left(\left|P^1_0\right|\right) +
t\, g\left(|\tilde P|\right) \right\}\, \mathrm{d}p &\le& (1+t)\, 
\int_{\mathbb{R}} |p|^k \max_i g_{t\,\left\|F_i\right\|_t}(p)\, 
\mathrm{d}p\\ &\le& (1+t)\, \max_i 
\int_{\mathbb{R}} |p|^k g_{t\,\left\|F_i\right\|_t}(p)\, \mathrm{d}p.
\end{eqnarray*}
Inequality \eqref{eq:mk1-mk2} is obtained by using the last lines
in Lemma~\ref{lem:moments}.
\end{proof}

A small variation on the above arguments allows us to prove that 
in fact, the density and flux are regular enough to satisfy the
continuity equation. This is important in order to be able to say that
the description via the Poisson equation~\eqref{gauss.dim.2} is 
equivalent to the Amp\`ere equation~\eqref{amp.dim.2}.

\begin{corollary}
Under Hypotheses~\mbox{\rm\ref{hypo:f0c1}--\ref{hypo:f0}} and the 
assumption~\eqref{lip-reg}, resp.~\eqref{cb-reg},
$n$ and $j$ are Lipschitz, resp.~continuously differentiable in space 
and time and they satisfy
\begin{equation}
\frac{\partial n}{\partial t} + \frac{\partial j}{\partial x} = 0 .
\label{eq:cons:mass}
\end{equation}
Moreover, there is a constant~$C$, depending only on the majorising 
function
$g(p)$, such that:
\begin{equation}
\left\|\frac{\partial n}{\partial x}\right\|_t \le C\, 
\mathrm{e}^{t\, \left(1+\left\|\partial_x F\right\|_t\right) }.
\label{46.4}
\end{equation}
\label{cor:rg.n+j}
\end{corollary}
\begin{proof}
When the solution is classical, we can apply the chain rule 
and estimate directly the derivative in $x$ of the solution $f(t,x,p)$ 
to get
$$
\left| \frac{\partial f}{\partial x} \right|\, \le \,
\left(\sup_x \left| \frac{\partial X}{\partial x}(0) \right|\right)\, 
\left| \frac{\partial f_0}{\partial x_0} \right| + 
\left( \sup_x \left| \frac{\partial P}{\partial x}(0)\right|\right)\, 
\left| \frac{\partial f_0}{\partial p_0} \right| .
$$
Lemma \ref{lem:char:der} and Hypothesis~\ref{hypo:f0} imply that
(see previous Lemma):
$$\left| \frac{\partial f}{\partial x} \right|\, \le \, 
\mathrm{e}^{t\, \left(1+\left\|\partial_x F\right\|_t\right) }\, 
g_{t\,\|F\|_t}(p),$$
and thus $\partial_x f$ is integrable. Moreover, 
Lebesgue's dominated convergence theorem implies that $n$ is 
Lipschitz or differentiable with respect to $x$ and hence~(\ref{46.4}).
Similar arguments apply to $\partial_t n$ and $\partial_x j$ and the 
continuity 
equation~(\ref{eq:cons:mass}) becomes an easy consequence of 
Eq.~\eqref{vla.adim} upon integration on $p$ and using 
Hypotheses~\ref{hypo:f0c1}--\ref{hypo:f0} on the initial data.
In case the solution is only mild under assumption~\eqref{lip-reg}, 
one can reproduce the Lipschitz bounds by estimating
the difference $n(t,x_1)-n(t,x_2)$ using analogous arguments to 
previous Lemma \ref{lem:difmoments}, we leave the details to the 
reader.
\end{proof}

And it stems from Lemma~\ref{lem:perio:X.P} that:

\begin{lemma}
Under the hypotheses:
$$\forall (x,p),\ f_0(x+L,p) = f_0(x,p)\,;\quad
\forall (t,x),\ \{E,A\}(t,x+L)=\{E,A\}(t,x)\,;$$
there holds:
$$\forall (t,x,p),\ f(t,x+L,p) = f(t,x,p),\ 
\{n,j\}(t,x+L)=\{n,j\}(t,x).$$
\label{lem:perio:n.j}
\end{lemma}

Let us finally remark that thanks to the mass conservation property 
of the Vlasov equation, we estimate the integrals of~$n$. 
In the open-space case, we have
\begin{equation}
\int_a^b n(x)\,\mathrm{d}x \le M,
\label{eq:int.n.open}
\end{equation}
where $M$~is the total mass of~$f_0$. In the periodic case,
\begin{equation}
\int_a^b n(x)\,\mathrm{d}x \le M\,\left\lceil\frac{b-a}{L}\right\rceil,
\label{eq:int.n.perio}
\end{equation}
where $M$~is now the mass of~$f_0$ over one period, and 
$\lceil r\rceil$ is the smallest integer larger or equal to~$r$.

Moreover, in the open-space case we can estimate moments in $x$ under
suitable additional assumptions on the initial data. In fact, let us 
consider the following
\begin{hypothesis}
$f_0(x,p)\le g(x)\, g(p)$ in ${\mathbb{R}}^2$.
\label{hypo:f0open}
\end{hypothesis}
By following the same lines of argument as in Lemma~\ref{lem:moments},
we can prove:
\begin{lemma}
Let $f$ be a characteristic solution of the Vlasov 
equation~\eqref{vla.adim} in 
the open-space case with force field $F \in \cesp{0}{\lip{1}}$ and 
assume that
Hypothesis~\mbox{\rm\ref{hypo:f0open}} is satisfied. 
Then, there exists a constant $Z$,
depending polynomially on the first two moments of $g(p)$, $t$ and 
$\left\|F \right\|_t$, such that
\begin{equation}
\int_{{\mathbb{R}}^2} (|x|+|p|)\,f(t,x,p)\, \mathrm{d}x\,\mathrm{d}p 
\le Z\left(M_0,M_1,t,\left\|F \right\|_t\right).
\label{lossofmass}
\end{equation}
\label{lem:momentsx}
\end{lemma}
This can be generalised to all moments of $f$ and its derivatives, 
under
\begin{hypothesis}
$\left|\partial_x f_0(x,p)\right| \le g(x)\,g(p)$ and 
$\left|\partial_p f_0(x,p)\right| \le g(x)\,g(p)$ in~${\mathbb{R}}^2$.
\label{hypo:f0open:bis}
\end{hypothesis}

\begin{lemma}
Assume that Hypotheses~{\rm\ref{hypo:f0open}} 
and~{\rm\ref{hypo:f0open:bis}} hold, and that $g$ has its moment of
order~$m$ bounded. Then, for all $k,\ \ell\le m$, there exists a 
constant $Z_{k,\ell}$ such that, for all $t\in[0,T]$:
\begin{eqnarray*}
\int\!\!\!\int \left(|x|^k + |p|^\ell\right)\, f(t,x,p)\, 
\mathrm{d}x\, \mathrm{d}p &\le& Z_{k,\ell},\\
\int\!\!\!\int \left(|x|^k + |p|^\ell\right)\, \left|\partial_x 
f(t,x,p) \right|\, \mathrm{d}x\, \mathrm{d}p &\le& Z_{k,\ell}\,
\mathrm{e}^{t\,\left(1+\|\partial_x F\|_t\right)},
\end{eqnarray*}
and a similar bound holds for $\partial_p f$.
\label{lem:momentsx:bis}
\end{lemma}


%

\section{Iterative procedure and global weak solutions%
\label{sec:fixptweak}}

In this section, we present an iterative procedure to solve the
1D Vlasov--Maxwell system~(\ref{vla.adim}--\ref{ond.adim}) for the NR 
and QR~cases. 

First, we define the iterative procedure based on the 
Cooper--Klimas\cite{CoKl80} approach. Then, we will derive estimates 
on the fields that allow us to obtain a limit by telescopic series. 
However, these estimates will not allow us to get a global 
characteristic solution, due to the lack of a global-in-time estimate 
on the space derivative of the density, which is needed to control 
that of the right-hand side in~\eqref{ond.adim}, and thus the second 
space derivative of the vector potential. Therefore, at this level of 
generality we are only able to obtain global weak solutions.
Improvements of this basic result, namely local (in the NR~case) and
global (in the QR~case) existence of classical and mild solutions
will be postponed to the next two sections.

\medbreak

We now fix an initial condition $f_0$ for the distribution function,
as well as initial data $(A_0,{\dot A}_0)$ for the vector 
potential. $n_0$ and~$E_0$ are the density and electrostatic field 
given by~$f_0$, as defined in~\eqref{ci2.adim}. Let us remind that we 
always assume that the external density verifies $n_\mathrm{ext}\in 
L^\infty({\mathbb{R}})$.
We also fix a target time~$T$, and, if we are interested in periodic
solutions, a space period~$L$. 

\subsection{Definition of the recurrence operator\label{ssc:df.rec}}
Given the field pair $(E,A)\in\cesp{0}{\lip{1}\times \lip{2}}$,
one constructs $(E',A')={\cal L}(E,A)$ as follows.
\begin{enumerate}
\item The characteristic system~(\ref{carac}) is solved, with the
force $F$ given by~(\ref{df:force}).
\item One computes the characteristic solution to Vlasov's
equation by~(\ref{eq:df:char}) and its density and flux~$(n,j)$
by~(\ref{df:dens:cour}).
\item Finally, $E'$ and $A'$ are computed as:
\begin{eqnarray}
E'(t,x) &:=& E_0(x) + \int_0^t j(s,x)\, \mathrm{d}s,
\label{43.1}\\
A'(t,x) &:=& \frac12\,\biggl\{ A_0(x+t) + A_0(x-t) + 
\int_{x-t}^{x+t} {\dot A}_0(y)\,\mathrm{d}y \nonumber\\
&&\hphantom{\frac12\,\biggl\lbrace}\mbox{}
-\int_{0}^t\int_{x+s-t}^{x+t-s} (n\,A)(s,y)\, \mathrm{d}y\,\mathrm{d}s 
\biggr\}.
\label{43.2}
\end{eqnarray}
\end{enumerate}
{F}rom Corollary~\ref{cor:rg.n+j} and the Duhamel 
formulae~(\ref{40.3}--\ref{40.5}), we immediately deduce:

\begin{theorem}
If $(A_0,{\dot A}_0)\in\lip{2}\times\lip{1}$, the 
operator~${\cal L}$ maps $\cesp{0}{\lip{1}\times \lip{2}}$ to itself.
If $(A_0,{\dot A}_0)\in\cb{2}\times\cb{1}$, 
$f_0\in{\cal C}^1({\mathbb{R}}^2)$, 
and $n_\mathrm{ext}\in\cb{0}$, then ${\cal L}$
maps $\cesp{0}{\cb{1}\times \cb{2}}$ to itself.

Moreover, the Poisson equation \eqref{poi.adim} is satisfied for 
the pair $(E',n)$; and the time derivative of~$A'$ is bounded in $x$, 
\hbox{i.e.}~$A'\in \cesp{1}{L^\infty({\mathbb{R}})}$ or 
$\cesp{1}{\cb{0}}$, even if $A$ does not belong \emph{a priori} to 
such a space.
\end{theorem}

This ensures that the operator ${\cal L}$ can be iterated.
Of course, we shall need some quantitative estimates. To 
establish them is the goal of the next Subsection.

Moreover, we deduce from Lemma~\ref{lem:perio:n.j} the following

\begin{corollary}
If $(A_0,{\dot A}_0)\in\lipp{2}\times\lipp{1}$, the 
operator~${\cal L}$ maps $\cesp{0}{\lipp{1}\times \lipp{2}}$ to itself.
If $(A_0,{\dot A}_0)\in\cbp{2}\times\cbp{1}$, 
$f_0\in{\cal C}^1_L({\mathbb{R}}^2)$, and $n_\mathrm{ext}\in\cbp{0}$,
then ${\cal L}$ maps $\cesp{0}{\cbp{1}\times \cbp{2}}$ to itself.
\label{cor:sol:per}
\end{corollary}

\subsection{\emph{A priori} estimates\label{ssc:apriori}}
In the sequel, we shall always assume \emph{at least} that
$(A_0,{\dot A}_0)\in\lip{2}\times\lip{1}$. The constants denoted 
$C$ or $C_i$ may vary from one line to the next, and depend on the 
initial conditions, $T$ and~$L$ (but on nothing else).

First, we fix $(E,A)\in\cesp{0}{\lip{1}\times \lip{2}}$, and
$(E',A'):={\cal L}(E,A)$. $F$ is the force field corresponding 
to~$(E,A)$.

For the electrostatic field we have the following properties:

\begin{lemma}
The following estimates hold:
\begin{eqnarray}
\left\|E'\right\|_t &\le& C_0 + C_1\, \int_0^t \|F\|_s\, \mathrm{d}s\,;
\label{44.1b} \\
\left\|\frac{\partial E'}{\partial x}\right\|_t &\le& C_0 + C_1\, 
\|F\|_t .
\label{44.2}
\end{eqnarray}
\label{lem:est:E'}
\end{lemma}
\begin{proof}
To obtain the first estimate, we proceed as in \cite{Bostan} by 
duality. Given
$\varphi \in L^1({\mathbb{R}})$, we apply changes of variables to get
\begin{eqnarray*}
\int_{\mathbb{R}} \int_0^t j(s,x)\, \varphi(x) \, \mathrm{d}s 
\,\mathrm{d}x & = & 
\int_{{\mathbb{R}}^2} f_0(x,p) \, \int_0^t P(s;0,x,p)\, 
\varphi(X(s;0,x,p)) \, \mathrm{d}s \,\mathrm{d}x \, \mathrm{d}p \\
&=& \int_{{\mathbb{R}}^2} f_0(x,p)\, \int_x^{X(t;0,x,p)} \varphi(u) \, 
\mathrm{d}u \,\mathrm{d}x \, \mathrm{d}p ,
\end{eqnarray*}
and thus, using Hypothesis \ref{hypo:f0}
\begin{eqnarray*}
\left| \int_{\mathbb{R}} \int_0^t j(s,x) \varphi(x) \, \mathrm{d}s 
\,\mathrm{d}x \right| 
&\le& \int_{{\mathbb{R}}^2} f_0(x,p) \, \left| \int_x^{X(t;0,x,p)} 
\varphi(u) \, \mathrm{d}u \right| \,\mathrm{d}x \, \mathrm{d}p \\
&\le& \int_{\mathbb{R}} g(p) \, \int_{\mathbb{R}} |\varphi(u)| \, 
\int_{\mathbb{R}} \chi(u;S(t,x,p))\, \mathrm{d}x \,
\mathrm{d}u \,\mathrm{d}p
\end{eqnarray*}
where $\chi(u;S(t,x,p))$ is the characteristic function of the set 
$S(t,x,p)=\{|u-x|\leq |X(t;0,x,p)-x|\}$ as a function of $u$.
Using the bound on the divergence of forward-in-time characteristics 
as in Lemma \ref{42.2}, we obtain 
$$ |X(t;0,x,p)-x|\leq t \int_0^t \|F\|_s \, \mathrm{d}s
\quad
\mbox{and}
\quad
\int_{\mathbb{R}} \chi(u;S(t,x,p)) \, \mathrm{d}x \leq 2 t\, \int_0^t 
\|F\|_s \, \mathrm{d}s\,;$$
and therefore, we deduce
\[
\left\| \int_0^t j(s,x) \,\mathrm{d}s \right\|_t \leq 2 t\, 
\int_{\mathbb{R}} g(p) \, \mathrm{d}p \int_0^t \|F\|_s \, \mathrm{d}s.
\label{est:j}
\]
Finally, \eqref{44.1b} is deduced directly from \eqref{43.1}.

Now, Corollary \ref{cor:rg.n+j} assures that the continuity 
equation~\eqref{eq:cons:mass} is satisfied. Taking into account that
equation and the Amp\`ere equation~\eqref{43.1}, then 
$$
\frac{\partial E'}{\partial x} = n_\mathrm{ext} - n,
$$ 
holds. Estimate \eqref{44.2} follows from \eqref{44.3} and 
$n_\mathrm{ext}\in L^\infty({\mathbb{R}})$.
\end{proof}

For the vector potential, we have a similar result.

\begin{lemma}
The following estimates hold:
\begin{eqnarray}
\left\|A'\right\|_t &\le& C_0 + C_1\, \int_0^t \|A\|_s\, \mathrm{d}s,
\label{57.1}\\
\left\|\frac{\partial A'}{\partial x}\right\|_t &\le& C_0 + \|A\|_t\,
\left( C_1 + C_2\, \int_0^t \|F\|_s\, \mathrm{d}s \right),
\label{50.2}\\
\left\|\frac{\partial A'}{\partial t}\right\|_t &\le& C_0 + \|A\|_t\,
\left( C_1 + C_2\, \int_0^t \|F\|_s\, \mathrm{d}s \right).
\label{50.2.bis}
\end{eqnarray}
\label{lem:est:A'}
\end{lemma}
\begin{proof}
{F}rom~(\ref{43.2}), we bound
$$|A'(t,x)|\le C_0 + \int_{0}^t \|A\|_s\, \left( 
\int_{x+s-t}^{x+t-s} n(s,y)\, \mathrm{d}y\right)\,\mathrm{d}s.
$$
Using (\ref{eq:int.n.open}) or~(\ref{eq:int.n.perio}), we then
bound the integral of~$n$ by $M$ or $M\,\lceil 2\,T/L \rceil$.
Hence~(\ref{57.1}).

Let us estimate the derivatives of $A$. $\partial_x A'$ is given by 
the Duhamel formula~(\ref{40.4}), with $f$ replaced by $-n\,A$. 
Hence the majoration:
$$
2\,\left|\frac{\partial A'}{\partial x}(t,x)\right| \le C_0 + \|A\|_t\,
\int_0^t (n(s,x+s-t)+n(s,x-s+t))\, \mathrm{d}s.
$$
Now, using the uniform bound on the density \eqref{44.3}, we obtain
$$\int_0^t (n(s,x+s-t)+n(s,x-s+t))\,\mathrm{d}s \le C_1 + C_2\, 
\int_0^t \|F\|_s\,\mathrm{d}s.$$
This gives~(\ref{50.2}).
As the Duhamel formula~(\ref{40.3}) for the time derivative is
very similar to~(\ref{40.4}), one can establish~(\ref{50.2.bis})  
by the same reasoning.
\end{proof}

Now, we consider two field pairs $(E_1,A_1)$ and $(E_2,A_2)$, 
with the corresponding forces $F_1$ and~$F_2$, and we set
$(E'_1,A'_1):={\cal L}(E_1,A_1)$, $(E'_2,A'_2):=
{\cal L}(E_2,A_2)$.

\begin{lemma}
The following estimates hold for all field pairs $(E_1,A_1)$ 
and $(E_2,A_2)$:
\begin{eqnarray}
\left\|E'_1 - E'_2\right\|_t &\le& R_1\left(M_0+M_1, 
\max_i \left\|F_i\right\|_t \right)\,
\int_0^t \|F_1-F_2\|_s\, \mathrm{d}s,
\label{44.1}\\
\left\|A'_1 - A'_2\right\|_t &\le& \left\|A_2\right\|_t\,
\left( C_0 + C_1\,\max_i \left\|F_i\right\|_t \right)\,
\int_0^t \left\|F_1 - F_2\right\|_s\, \mathrm{d}s 
\nonumber\\
&&\mbox{} + C_2\,\int_{0}^t \left\|A_2-A_1\right\|_s\, \mathrm{d}s.
\label{55.4}
\end{eqnarray}
\label{lem:47.3}
\end{lemma}
\vspace{-12pt}
\begin{proof}
{F}rom \eqref{43.1} we deduce
$$\left\|E'_1 - E'_2\right\|_t \le \int_0^t |j_1(s,x)-j_2(s,x)| \, 
\mathrm{d}s,$$
and thus \eqref{44.1} is a simple consequence of 
Lemma~\ref{lem:difmoments}.

{F}rom~(\ref{43.2}), we have
$$2\,(A'_1-A'_2)(t,x) = \int_0^t \int_{x+s-t}^{x+t-s} 
(n_2\,A_2 - n_1\,A_1)(s,y)\, \mathrm{d}y\,\mathrm{d}s.$$
Writing: $n_2\,A_2 - n_1\,A_1 = (n_2-n_1)\,A_2 - (A_2 - A_1)\,n_1$,
we bound:
\begin{eqnarray*}
2\,\left\|A'_1 - A'_2\right\|_t &\le&
\int_0^t \left\|A_2\right\|_s\, \int_{x+s-t}^{x+t-s} 
\left| n_2(s,y)-n_1(s,y) \right|\, \mathrm{d}y\, \mathrm{d}s 
\nonumber\\
&&\mbox{}+ \int_0^t 
\left\|A_2-A_1\right\|_s\, \int_{x+s-t}^{x+t-s} 
n_1(s,y)\, \mathrm{d}y\, \mathrm{d}s. 
\end{eqnarray*}
Using (\ref{eq:int.n.open}) or~(\ref{eq:int.n.perio}),
the second line of this inequality is easily bounded as
$\lambda\,M\,\int_{0}^t \left\|A_2-A_1\right\|_s\, \mathrm{d}s$, with 
$\lambda=1$ or $\lceil2\,T/L\rceil$.
Then, the first line is bounded thanks to Lemma~\ref{lem:difmoments}:
\begin{equation}
\left|n_1 - n_2\right|(t,x) \le \left( C_0 + 
C_1\,\max_i \left\|F_i\right\|_t \right)\,
\int_0^t \left\|F_1 - F_2\right\|_s\, \mathrm{d}s,
\label{eq:n1-n2:2}
\end{equation}
and finally:
\begin{eqnarray*}
\left\|A'_1 - A'_2\right\|_t &\le& \left\|A_2\right\|_t\,
\left( C_0 + C_1\,\max_i \left\|F_i\right\|_t \right)\,
\int_0^t \left\|F_1 - F_2\right\|_s\, \mathrm{d}s
\\
&&\mbox{} + \frac{\lambda\,M}2\,\int_{0}^t 
\left\|A_2-A_1\right\|_s\, \mathrm{d}s,
\end{eqnarray*}
which is~(\ref{55.4}).
\end{proof}

Now, we estimate the difference of the derivatives.

\begin{lemma}
There holds, for any $(E_1,A_1)$ and $(E_2,A_2)$:
\begin{eqnarray}
\left\|\frac{\partial E'_1}{\partial x} - 
\frac{\partial E'_2}{\partial x}\right\|_t 
&\le& \left( C_0 + C_1\,\max_i \left\|F_i\right\|_t \right)\,
\int_{0}^t \left\|F_2-F_1\right\|_s\, \mathrm{d}s ,
\label{52.1e}\\
\left\|\frac{\partial A'_1}{\partial x} - 
\frac{\partial A'_2}{\partial x}\right\|_t 
&\le& \left( C_0 + C_1\,\max_i \left\|F_i\right\|_t \right)\,
\biggl\{ \int_{0}^t \left\|A_2-A_1\right\|_s\, \mathrm{d}s \nonumber\\
&&\hspace*{8em}\mbox{} + \left\|A_2\right\|_t\, 
\int_0^t \left\|F_1 - F_2\right\|_s\, \mathrm{d}s \biggr\},
\label{52.1}\\
\left\|\frac{\partial A'_1}{\partial t} - 
\frac{\partial A'_2}{\partial t}\right\|_t 
&\le& \left( C_0 + C_1\,\max_i \left\|F_i\right\|_t \right)\,
\biggl\{ \int_{0}^t \left\|A_2-A_1\right\|_s\, \mathrm{d}s \nonumber\\
&&\hspace*{8em}\mbox{} +
\left\|A_2\right\|_t\, \int_0^t \left\|F_1 - F_2\right\|_s\, 
\mathrm{d}s \biggr\}.
\label{52.1.bis}
\end{eqnarray}
\label{lem:51.1}
\end{lemma}
\vspace{-12pt}
\begin{proof}
The Poisson equation and (\ref{eq:n1-n2:2}) imply (\ref{52.1e}).
{F}rom the Duhamel formula~(\ref{40.4}), with $f$~replaced
successively with $-n_1\,A_1$ and~$-n_2\,A_2$, we derive:
\begin{eqnarray*}
2\,\left(\frac{\partial A'_1}{\partial x} - 
\frac{\partial A'_2}{\partial x}\right)(t,x) 
&=& \int_0^t (n_2\,A_2 - n_1\,A_1)(s,x+s-t)\,\mathrm{d}s \\
&&\mbox{} - \int_0^t (n_2\,A_2 - n_1\,A_1)(s,x+t-s)\,\mathrm{d}s \;;
\end{eqnarray*}
once more, we write
$n_2\,A_2 - n_1\,A_1 = (n_2-n_1)\,A_2 - (A_2 - A_1)\,n_1$,
which gives
\begin{eqnarray*}
\left\|\frac{\partial A'_1}{\partial x} - 
\frac{\partial A'_2}{\partial x}\right\|_t
&\le& \left\|A_2\right\|_t\, \int_0^t 
\left\|n_1 - n_2\right\|_s\, \mathrm{d}s +
\left\|n_1\right\|_t\, \int_0^t 
\left\|A_1 - A_2\right\|_s\, \mathrm{d}s.
\end{eqnarray*}
Using the bounds~(\ref{eq:n1-n2:2}) for the first term, 
and~(\ref{44.3}) for the second yields
\begin{eqnarray*}
\left\|\frac{\partial A'_1}{\partial x} - 
\frac{\partial A'_2}{\partial x}\right\|_t
&\le& \left\|A_2\right\|_t\,
\left( C_0 + C_1\,\max_i \left\|F_i\right\|_t \right)\,
\int_0^t \left\|F_1 - F_2\right\|_s\, \mathrm{d}s \\
&&\mbox{} + \left( C_2\, \|F_1\|_t + C_3 \right)\,\int_0^t 
\left\|A_1 - A_2\right\|_s\, \mathrm{d}s,
\end{eqnarray*}
which implies~(\ref{52.1}).
Once again, the similarity of the formulae (\ref{40.3}) 
and~(\ref{40.4}) allows to deduce~(\ref{52.1.bis}).
\end{proof}

\subsection{Convergence of successive approximations}
We start from the initial data $(E_0(x),A_0(x))$, by extending
them to constant-in-time functions over $[0,T]\times{\mathbb{R}}$.
Then, we construct a sequence
$(E_k,A_k)_{k\in{\mathbb{N}}}$ by the recurrence formula:
$$(E_{k+1},A_{k+1}):= {\cal L}(E_k,A_k),\quad \forall k\ge0.$$
Of course, we set $F_k:= E_k + A_k\,\partial_x A_k$.

To establish the convergence of $(E_k,A_k)_{k\in{\mathbb{N}}}$,
the following result will be useful. It is easily proved by
induction.
\begin{lemma}
Let $(u_k)_{k\in{\mathbb{N}}}$ be a sequence of positive functions,
$u_k:\ [0,T]\to{\mathbb{R}}^+$, satisfying
\begin{eqnarray*}
\mbox{\rm(i)}\quad&&\forall t\in[0,T],\quad u_0(t)\le c, \\
\mbox{\rm(ii)}\quad&& \forall k\in{\mathbb{N}},\ \forall t\in[0,T],
\quad u_{k+1}(t) \le a + b\,\int_0^t u_k(s) \, \mathrm{d}s,
\end{eqnarray*}
for some constants $a,\ b,\ c\in{\mathbb{R}}^+$.
Then the following estimate holds:
$$\forall k\in{\mathbb{N}},\ \forall t\in[0,T],\quad
u_k(t) \le a\,\sum_{i=1}^{k-1} \frac{b^i\,t^i}{i!} + 
c\,\frac{b^k\,t^k}{k!} \;;$$
hence the whole sequence is uniformly bounded by a constant $u_*$
which depends on $a$, $b$, $c$ and~$T$.
\label{lem:56.2}
\end{lemma}

\begin{theorem}
The sequence $(E_k,A_k)_{k\in{\mathbb{N}}}$ converges uniformly in $t$ 
and~$x$ towards a limit~$(E,A)$.
\begin{enumerate}
\item If $(A_0,{\dot A}_0)\in\lip{2}\times\lip{1}$, the sequences 
$(E_k)$ and~$(A_k)$ converge respectively within $\cesp{0}{\lip{1}}$
and $\cesp{1}{L^\infty({\mathbb{R}})}\cap \cesp{0}{\lip{1}}$.
\item If $(A_0,{\dot A}_0)\in\cb{2}\times\cb{1}$, $f_0\in
{\cal C}^1({\mathbb{R}}^2)$, 
and $n_\mathrm{ext}\in\cb{0}$, then the convergences take place
within $\cesp{0}{\cb{1}}$ and $\cesp{1}{\cb{0}}\cap \cesp{0}{\cb{1}}$.
\item The above conclusions are valid when all 
spaces are replaced by their periodic counterparts.
\end{enumerate}
\label{th:cv:fort}
\end{theorem}
\begin{proof}
First, Eq.~(\ref{57.1}) gives: 
$$\left\|A_{k+1}\right\|_t\le C_0 +
C_1\,\int_0^t \left\|A_{k}\right\|_s\, \mathrm{d}s$$
and thus, the previous lemma ensures that  
$\left\|A_k\right\|_t\le A_*$ uniformly in $k$ and~$t$.
Then, Eqs. (\ref{44.1b}) and~(\ref{50.2}) imply:
$$\left\|F_{k+1}\right\|_t\le \left\|E_{k+1}\right\|_t + A_*\,
\left\|\partial_x A_{k+1}\right\|_t 
\le C_2 + C_3\, \int_0^t \left\|F_k\right\|_s\, \mathrm{d}s,$$
so, once more, $\left\|F_k\right\|_t\le F_*$. 
Using again~(\ref{44.1b}) and~(\ref{50.2}) shows that the sequences 
$\left\|E_k\right\|_t$ and $\left\|\partial_x A_k\right\|_t$ are 
uniformly bounded in $k$ and~$t$. 

All these bounds allow to rewrite~(\ref{55.4}) and~(\ref{52.1}) as
\begin{eqnarray}
\left\|A_{k+1} - A_k\right\|_t &\le& 
C_1\, \int_0^t \left\|F_k - F_{k-1}\right\|_s\, \mathrm{d}s + 
C_2\, \int_{0}^t \left\|A_k-A_{k-1}\right\|_s\, \mathrm{d}s \;;\;
\label{bnd:Ak:sq}\\
\left\|\frac{\partial A_{k+1}}{\partial x} - 
\frac{\partial A_k}{\partial x}\right\|_t 
&\le& C_1\,\int_0^t \left\|F_k - F_{k-1}\right\|_s\, \mathrm{d}s +
C_2\, \int_{0}^t \left\|A_k-A_{k-1}\right\|_s\, \mathrm{d}s .\qquad
\label{bnd:Bk:sq}
\end{eqnarray}
Then, writing
$$F_{k+1} - F_k = E_{k+1} - E_k + A_{k+1}\,
\left(\partial_x A_{k+1}-{\partial_x A_k}\right) + {\partial_x A_k}\,
\left(A_{k+1} - A_k\right),$$
one deduces from~(\ref{44.1}), (\ref{bnd:Ak:sq}) and~(\ref{bnd:Bk:sq})
\begin{equation}
\left\|F_{k+1} - F_k\right\|_t \le 
C_1\, \int_0^t \left\|F_k - F_{k-1}\right\|_s\, \mathrm{d}s + 
C_2\, \int_{0}^t \left\|A_k-A_{k-1}\right\|_s\, \mathrm{d}s.
\label{bnd:Fk:sq}
\end{equation}
Hence, the sequence $u_k(t):= \left\|A_k - A_{k-1}\right\|_t + 
\left\|F_k - F_{k-1}\right\|_t$ satisfies: 
$u_{k+1}(t) \le b\,\int_0^t u_k(s) \, \mathrm{d}s$, for some 
constant~$b$. Lemma~\ref{lem:56.2} with $a=0$ then implies 
$u_k(t) \le b^k\,T^k/k!$, \hbox{i.e.}:
$$\left\|A_{k+1} - A_k\right\|_t \le\frac{b^k\,T^k}{k!},\quad
\left\|F_{k+1} - F_k\right\|_t\le\frac{b^k\,T^k}{k!},$$
so by~(\ref{44.1}) and~(\ref{bnd:Bk:sq})
$$\left\|E_{k+1} - E_k\right\|_t \le C\,\frac{b^k\,T^k}{k!},\quad
\left\|\frac{\partial A_{k+1}}{\partial x} - 
\frac{\partial A_k}{\partial x}\right\|_t 
\le C\,\frac{b^k\,T^k}{k!}.$$
As a consequence, the four sequences 
$\left(E_k\right)_{k\in{\mathbb{N}}}$, 
$\left(A_k\right)_{k\in{\mathbb{N}}}$, 
$\left(\partial_x A_k\right)_{k\in{\mathbb{N}}}$, 
$\left(F_k\right)_{k\in{\mathbb{N}}}$, all converge uniformly on
$[0,T]\times{\mathbb{R}}$. Let $E$, $A$, $B$, $F$, be the limits; 
clearly $B=\partial_x A$ in the sense of distributions, and $F = E + 
A\,\partial_x A$.  Under the hypothesis 
$(A_0,{\dot A}_0)\in\lip{2}\times\lip{1}$, each term of the four sequences 
is in~$\cesp{0}{L^\infty({\mathbb{R}})}$, and so are the limits.
So: $E\in\cesp{0}{L^\infty({\mathbb{R}})}$ and $A\in\cesp{0}{\lip{1}}$.
If $(A_0,{\dot A}_0)\in\cb{2}\times\cb{1}$, $f_0\in{\cal C}^1({\mathbb{R}}^2)$, 

and $n_\mathrm{ext}\in\cb{0}$, the terms and the limits are
in~$\cesp{0}{\cb{0}}$, \hbox{i.e.}~$E\in\cesp{0}{\cb{0}}$ and 
$A\in\cesp{0}{\cb{1}}$.

With all these results, Eq.~(\ref{52.1.bis}) shows that:
$\left\|\partial_t A_{k+1}-{\partial_t A_k}\right\| \le 
C\,b^k\,T^k/k!$; since $A_0$ is independent of time, the sequence 
$\left(\partial_t A_k\right)_{k\in{\mathbb{N}}}$ also converges 
uniformly in $\cesp{0}{L^\infty({\mathbb{R}})}$, 
resp.~$\cesp{0}{\cb{0}}$, towards a limit
which is necessarily equal to $\partial_t A$. Thus, $A\in
\cesp{1}{L^\infty({\mathbb{R}})}$, resp.~$\cesp{1}{\cb{0}}$.

Similarly, Eq.~(\ref{52.1e}) implies that 
$\left(\partial_x E_k\right)_{k\in{\mathbb{N}}}$ 
is a Cauchy sequence and thus converges uniformly toward a limit which 
is necessarily equal to~$\partial_x E$. Hence, $E\in
\cesp{0}{\lip{1}}$ or $\cesp{0}{\cb{1}}$.

The last point easily follows from Corollary~\ref{cor:sol:per}.
\end{proof}

\subsection{Global existence and properties of weak solutions}
If $(E,A)\in\cesp{0}{\lip{1}}\times\cesp{0}{\lip{2}}$, then 
it is a fixed point of ${\cal L}$ within this space. Thus, the triple
$(f,E,A)$, where the function $f$ is defined by~(\ref{eq:df:char}),
is a characteristic solution to~(\ref{vla.adim}--\ref{poi.adim}). 

This fixed point and its associated characteristic solution, 
if they exist, are unique. Indeed, let 
$(f_1,E_1,A_1)$ and~$(f_2,E_2,A_2)$ be two such solutions, 
with respective force fields $F_1$ and~$F_2$.
Reasoning as in the proof of the above theorem shows:
$$\left\| A_1 - A_2 \right\|_t + \left\| F_1 - F_2 \right\|_t
\le C\,\int_0^t \left\{ \left\| A_1 - A_2 \right\|_s + 
\left\| F_1 - F_2 \right\|_s\right\}\,\mathrm{d}s\,;$$
hence $A_1=A_2$ and $F_1=F_2$ by Gronwall's lemma, then $E_1=E_2$ 
by~(\ref{44.1}), and finally $f_1=f_2$ by the Cauchy--Lipschitz
theorem.

Without further \emph{a priori} estimates on the second derivatives
of the vector potential $A$, we will not be able to obtain global
existence of characteristic solutions. However, the bounds in 
Subsection~\ref{ssc:apriori} show that the norm of ${\cal L}(E,A)$ 
in $\cesp{0}{L^\infty({\mathbb{R}})\times\lip{1}}$, resp. 
$\cesp{0}{\cb{0}\times\cb{1}}$, is controlled by the norm of $(E,A)$ 
in the same space. These considerations suggest to extend
the operator~${\cal L}$ to a bigger space in order to obtain
weak solutions, which would also enjoy a 
\emph{uniqueness property}. 
This programme cannot be achieved within the framework
of $\cesp{0}{L^\infty({\mathbb{R}})\times\lip{1}}$. The reason is the 
well-known lack of density of smooth functions within 
$L^\infty({\mathbb{R}})$ or~$\lip{1}$. But this obstruction is 
removed if one works within $\cesp{0}{\cb{0}\times\cb{1}}$.

\begin{lemma}
Assume $(A_0,{\dot A}_0)\in\cb{2}\times\cb{1}$, 
$f_0\in{\cal C}^1({\mathbb{R}}^2)$, and $n_\mathrm{ext}\in\cb{0}$.
The operator ${\cal L}$ can be extended to a continous operator from
$\cesp{0}{\cb{0}\times\cb{1}}$ to itself, which satisfies the
estimates of Subsection~\ref{ssc:apriori}.
If the initial conditions are periodic, then 
${\cal L}$ maps $\cesp{0}{\cbp{0}\times\cbp{1}}$ \mbox{to itself.}
\label{lem:68.1}
\end{lemma}
\begin{proof}
${\cal L}$ is defined on the dense subspace $Y:=
\cesp{0}{\cb{1}\times\cb{2}}$ of $X:=\cesp{0}{\cb{0}\times\cb{1}}$,
with values in~$X$. {F}rom 
(\ref{44.1},\ \ref{55.4},\ \ref{52.1}), we see that it is uniformly 
continuous, in the norm of~$X$, on any set $K\cap Y$, where $K$ is a 
bounded set of~$X$.
Hence, it admits a unique continuous extension from~$X$ to itself. 

The extension procedure preserves the estimates of 
Subsection~\ref{ssc:apriori}, as their r.h.s. are uniformly continuous 
in the norm of~$X$ on any $K\cap Y$. 
Finally, the closedness of $\cesp{0}{\cbp{0}\times\cbp{1}}$ within~$X$
guarantees the invariance of this space by the extended operator.
\end{proof}

Hence, the following theorem, whose proof rephrases that
of Theorem~\ref{th:cv:fort}.

\begin{theorem}
Assume $(A_0,{\dot A}_0)\in\cb{2}\times\cb{1}$, 
$f_0\in{\cal C}^1({\mathbb{R}}^2)$, 
and $n_\mathrm{ext}\in\cb{0}$.
The operator ${\cal L}$ admits a unique fixed point $(E,A)\in
\cesp{0}{\cb{0}\times\cb{1}}$, which moreover belongs to
$\cesp{0}{\cb{1}}\times\cesp{1}{\cb{0}}$. If the initial
conditions are periodic, so is the fixed point.
\label{th:cv:faible:1}
\end{theorem}

Let us now check that this unique fixed point defines a distribution 
function $f$ in such a way that the triple $(f,E,A)$ is a global
weak solution of the Vlasov--Maxwell 
system~(\ref{vla.adim}--\ref{poi.adim}).

\begin{theorem}
Under the hypotheses of Theorem~\ref{th:cv:faible:1},
let $(E_k,A_k)_k$ a sequence converging to the fixed point $(E,A)$; 
without loss of generality, we can assume that its terms belong 
to~$\cesp{0}{\cb{1}\times\cb{2}}$. Let $(f_k)_k$ be the associated 
sequence of distribution functions obtained by the method of
characteristics. Then, $(f_k)_k$ converges uniformly in all its 
variables toward a function $f\in
\cesp{0}{{\cal C}^0_b({\mathbb{R}}^2)}$, which 
does not depend on the sequence $(E_k,A_k)_k$, and the triple~$(f,E,A)$
satisfies~\eqref{vla.adim}--\eqref{poi.adim} in the sense of 
distributions.
\label{th:cv:faible:2}
\end{theorem}
\begin{proof}
Repeating the argument of Lemma~\ref{lem:difmoments},
we easily obtain the estimate
\begin{equation}
\left| f_k(t,x,p) - f_\ell(t,x,p) \right| \le (1+t)\,
\max_{i=k,\ell} g_{t\,\left\|F_i\right\|_t}(p)\,
\int_0^t \left\|F_k - F_\ell\right\|_s\, \mathrm{d}s,
\label{eq:cv:uni:fk}
\end{equation}
from which follows that $\left(f_k(t)\right)_k$ is a Cauchy sequence 
in $L^\infty$ norm, uniformly in~$t$.
A similar argument shows that, if we consider another approximating 
sequence denoted by tildes, $f_k(t) - \tilde f_k(t)$ will converge
to zero, uniformly in $x,\ p$ and~$t$. 
Let us then call $f$ the common limit of all the sequences
$\left(f_k\right)_k$; its continuity in all variables follows from 
that of the~$f_k$.

As $f_k$ is a classical solution to the Vlasov equation, it is also
a distributional solution, \hbox{i.e.}~an equation similar 
to~\eqref{eq:df:weak} 
holds, with $f$ and~$F$ replaced resp.~with $f_k$ and~$F_k$. 
For any test function~$\psi$, the integrals in this formula are
taken over a compact subset of $[0,T)\times{\mathbb{R}}^2$. Thus, the 
uniform convergence in $x,\ p$ and~$t$ of the sequences 
$\left(f_k\right)_k$ and $\left(F_k\right)_k$ ensures 
that~\eqref{eq:df:weak} will also
hold at the limit: the limiting triple $(f,E,A)$ 
satisfies~(\ref{vla.adim}--\ref{poi.adim}) in the sense of 
distributions.
\end{proof}

The above result can be improved in the following two ways. First, we
show that $f(t)$ will retain the mass of the initial 
distribution~$f_0$.
\begin{proposition}
In the periodic case, $f_k(t)$ converges toward~$f(t)$ 
in~$L^1_L({\mathbb{R}}^2)$, uniformly in~$t$, and the mass of $f(t)$ 
over one space period is equal to that of~$f_0$. 
In the open-space case, under Hypothesis~\ref{hypo:f0open},
$f_k(t)$ converges weakly toward~$f(t)$ 
in~$L^1({\mathbb{R}}^2)$, and the mass of $f(t)$ is equal to that 
of~$f_0$. 
\label{th:nolossofmass}
\end{proposition}
\begin{proof}
Inequality~\eqref{eq:cv:uni:fk} shows that $f_k(t,x,\cdot)$ converges 
to~$f(t,x,\cdot)$ in $L^1({\mathbb{R}})$, uniformly in $x$ and~$t$. 
Hence the convergence of $f_k(t)$ toward~$f(t)$ in 
$L^1((a,b)\times{\mathbb{R}})$ for any $a<b$. 

Yet, the $f_k$ are \emph{classical} solutions to the Vlasov 
equation. In the periodic case, they conserve the mass of~$f_0$ over 
one space period. So, $f(t)$ retains this mass.

In the open-space case, taking into account the $L^1$ convergence in 
bounded intervals in space, it suffices to show that mass does not 
escape at infinity, uniformly in~$t$. Using Lemma~\ref{lem:momentsx}, 
we obtain:
$$\int_{|x|+|p|>R} |f_k(t,x,p)|\,\mathrm{d}x\,\mathrm{d}p \le
\frac1R\, \int_{|x|+|p|>R} (|x|+|p|)\, f_k(t,x,p)\,\mathrm{d}x\,
\mathrm{d}p \le \frac{Z}R \to0,$$
uniformly in $k$ and~$t$ when $R\to\infty$. Thus, it is a standard 
argument to check that
$\int_{{\mathbb{R}}^2} f(t,x,p)\,\mathrm{d}x\,\mathrm{d}p = \lim_k 
\int_{{\mathbb{R}}^2} f_k(t,x,p)\,\mathrm{d}x\,\mathrm{d}p = M$. 
\end{proof}

The second useful precision is that $f(t)$ is ``almost a mild 
solution'' to~\eqref{vla.adim}.
\begin{theorem}
The characteristics given by the~$F_k$ converge uniformly in all their 
variables toward a limit~$(X,P)$, and Eq.~\eqref{eq:df:mild} holds for 
this~$(X,P)$ and $f(t)$ given by Theorem~\ref{th:cv:faible:2}.
\end{theorem}
\begin{proof}
Let $(X^k,P^k)$ be the characteristics associated to~$F_k$. Their
uniform convergence follows from~(\ref{42.2}--\ref{42.4}).
As the $f_k$ are classical solutions to the Vlasov equation, they are 
also mild solutions, \hbox{i.e.}~they satisfy
\begin{equation}
\int_{{\mathbb{R}}^2} \psi(x,p)\,f_k(t,x,p)\, \mathrm{d}x\,\mathrm{d}p 
= \int_{{\mathbb{R}}^2} \psi(X^k(t;0,x,p),P^k(t;0,x,p))\,f_0(x,p)\, 
\mathrm{d}x\,\mathrm{d}p
\label{eq:mild:k}
\end{equation}
for all test functions $\psi\in {\cal C}_c({\mathbb{R}}^2)$ and all 
$t>0$. 

The sequence $\psi(X^k(t;0,x,p),P^k(t;0,x,p))\,f_0(x,p)$, 
resp. $\psi(x,p)\,f_k(t,x,p)$, converge uniformly on the support of 
$\psi$ towards $\psi(X(t;0,x,p),P(t;0,x,p))\,f_0(x,p)$, 
resp.~$\psi(x,p)\,f(t,x,p)$.
Hence, the integrals on both sides of~\eqref{eq:mild:k} converge
toward the two sides of~\eqref{eq:df:mild}, and the latter equality is
obtained at the limit.
\end{proof}


\section{Local-in-time existence of characteristic solutions
\label{sec:classmall}}
In order to have stronger solutions, \hbox{i.e.}~such that $f$ 
satisfies~(\ref{vla.adim}) in the characteristic sense, we need 
estimates on the Lipschitz constant of the force, which in turn 
amounts to bounding the second derivative of the vector potential.
In this Section we show that, in both the NR and QR~cases,
the limiting vector potential given by Theorem~\ref{th:cv:fort} does
satisfy such a bound, at least for a short time.

In the sequel, the operator~${\cal L}$ will be that of
Subsection~\ref{ssc:df.rec}.
We no longer consider the extended version, so that we can obtain
estimates based on characteristics.

\begin{lemma}
Let $(E,A)\in\cesp{0}{\lip{1}\times\lip{2}}$, and
$(E',A')={\cal L}(E,A)$. Then, the second derivative of~$A'$ is 
bounded as:
\begin{equation}
\left\|\frac{\partial^2 A'}{\partial x^2}\right\|_t \le C_0 + t\, 
\left\{ \left\|\frac{\partial A}{\partial x}\right\|_t\,\left(C_1 + 
C_2\, \|F\|_t \right) + C_3\,\|A\|_t\, 
\mathrm{e}^{t\, \left(1+\left\|\partial_x F\right\|_t\right)} \right\}.
\label{52.3}
\end{equation}
\label{lem:52.2}
\end{lemma}
\begin{proof}
We use the formula~(\ref{40.5}), with $f$ replaced by 
$-n\,A$, and the bounds (\ref{44.3}) and~(\ref{46.4}) for~$n$:
\begin{eqnarray*}
2\,\left|\frac{\partial^2 A'}{\partial x^2}(t,x)\right| &\le& C_0 + 
\int_0^t \left(  n\,\left|\frac{\partial A}{\partial x}\right| 
+   |A|\,\left|\frac{\partial n}{\partial x}\right| 
\right)(s,x+s-t) \,\mathrm{d}s\\
&&\hphantom{C_0 +}\mbox{} +
\int_0^t \left(  n\,\left|\frac{\partial A}{\partial x}\right| 
+   |A|\,\left|\frac{\partial n}{\partial x}\right| 
\right)(s,x-s+t) \,\mathrm{d}s\\
&\le& C_0 + t\, \left\|\frac{\partial A}{\partial x}\right\|_t\,
\left(C_1 + C_2\, \|F\|_t \right) + t\, \|A\|_t\,C_3\,
\mathrm{e}^{t\, \left(1+\left\|\partial_x F\right\|_t\right)} .
\end{eqnarray*}
\end{proof}

\begin{theorem}
Let $(E_k,A_k)_{k\in{\mathbb{N}}}$ and $(E,A)$ be as in 
Theorem~\ref{th:cv:fort}. There exists $0<T^*\le T$ such that, for
$0\le t< T^*$, the sequence 
$\left(\partial_x^2 A_k\right)_{k\in{\mathbb{N}}}$ is uniformly 
bounded in $x$ and~$t$. 
As a consequence, the couple $(E,A)$ allows to define a characteristic 
solution~$f$ to~(\ref{vla.adim}) on the interval 
$\left\lbrack0,T^*\right)$.
\label{th:solcla:1}
\end{theorem}
\begin{proof}
We have:
$$ \partial_x F_k = \partial_x E_k + \left(\partial_x A_k\right)^2 + 
A_k\, \partial_x^2 A_k.$$
Thus, the boundedness results of Theorem~\ref{th:cv:fort} 
and~(\ref{52.3}) yield:
$$\left\|\partial_x F_{k+1}\right\|_t \le \alpha + \beta\, t\,
\mathrm{e}^{t\,\left(1+\left\|\partial_x F_{k}\right\|_t\right)},$$
for some positive constants $\alpha,\ \beta$ depending only on the 
initial conditions, $T$ and (possibly)~$L$. Hence, 
$\left\|\partial_x F_{k}\right\|_t \le v_k(t)$, where $v_k(t)$ is 
defined by the recurrence formula:
\begin{equation}
v_0(t) = \left\|\partial_x F_{0}\right\|_t,\quad
v_{k+1}(t) = \alpha + \beta\, t\,\mathrm{e}^{t\,\left(1+v_k(t)\right)} 
:= \varphi_t(v_k(t)).
\label{df:vnt}
\end{equation}
In~\ref{app:b} we show that, for $t<T^*$, the sequence
$(v_k(t))_{k\in{\mathbb{N}}}$ is convergent and hence bounded. This
gives the uniform boundedness of 
$(\partial_x F_k(t))_{k\in{\mathbb{N}}}$
and, by~(\ref{52.3}), that of 
$(\partial_x^2 A_k(t))_{k\in{\mathbb{N}}}$, on
$\left\lbrack0,T^*\right)\times{\mathbb{R}}$. 
Hence, the latter sequence admits a subsequence 
which converges weakly-$*$ in $L^\infty([0,T^*]\times{\mathbb{R}})$
towards a limit which is necessarily equal to $\partial_x^2 A$. 
This gives the last part of the conclusion.
\end{proof}

\section{Global characteristic solutions in the quasi-relativistic case
\label{sec:clasqr}}
In this section, \emph{we only consider the QR case}.
We show that, in this framework, characteristic solutions exist for
any time. This result rests on a subtler treatment of the 
formula~(\ref{40.5}). On the other hand, we need to introduce
another assumption: 
\begin{hypothesis}
The majorising function $g$ of Hypothesis~\ref{hypo:f0} moreover 
satisfies
$$\int |p|^2\,g(p)\,\mathrm{d}p < \infty.$$
\label{hypo:f0:bis}
\end{hypothesis}
With this hypothesis, we can derive a better bound on $\partial_x^2 A$.
\begin{lemma}
Let $E\in\cesp{0}{\lip{1}}$, $A\in\cesp{0}{\lip{1}}\cap
\cesp{1}{L^\infty({\mathbb{R}})}$, and
$(E',A')={\cal L}(E,A)$. Then, the second derivative of~$A'$ is 
bounded as:
\begin{equation}
\left\|\frac{\partial^2 A'}{\partial x^2}\right\|_t \le C_0 + 
C_1\,\Pi_3\left(t\,\|F\|_t\right)\,
\left\{ \left\|\frac{\partial A}{\partial x}\right\|_s + 
\left\|\frac{\partial A}{\partial t}\right\|_s + \|A\|_s\,(1+\|F\|_s) 
\right\},
\label{52.3:new}
\end{equation}
where $\Pi_3$ is a polynomial of third degree, whose coefficients are
positive and depend only on the function~$g$.
\label{lem:52.2:new}
\end{lemma}
\begin{proof}
We use the formula~(\ref{40.5}) with $f$ replaced by $-n\,A$: 
\begin{eqnarray}
2\,\frac{\partial^2 A'}{\partial x^2}(t,x) &\le& D_0(t,x) + 
\int_0^t\left\{ \left(n\,\frac{\partial A}{\partial x}\right)(s,x+s-t) 
- \left(n\,\frac{\partial A}{\partial x}\right)(s,x+t-s) \right\}\,
\mathrm{d}s \nonumber\\
&&\mbox{} - \int_0^t \left\{ \left(A\,\frac{\partial n}{\partial x}
\right)(s,x+s-t) + \left(A\,\frac{\partial n}{\partial x}
\right)(s,x+t-s) \right\}\,\mathrm{d}s
\label{eq:bnd:A''.1}
\end{eqnarray}
where $D_0$ depends on the initial data only. The first line is
bounded by~(\ref{44.3}) as:
$$C_0 + \left( C_1 + C_2\, t\,\|F\|_t \right)\,
\int_0^t \left\|\frac{\partial A}{\partial x}\right\|_s\,\mathrm{d}s.$$

In order to bound the second line in~(\ref{eq:bnd:A''.1}),
we use~(\ref{vla.adim}) to rewrite 
$$(\widehat p-1)\, \frac{\partial f}{\partial x}(s,x+s-t,p) = 
-\frac{\partial}{\partial s}[f(s,x+s-t,p)] + 
F(s,x+s-t)\,\frac{\partial f}{\partial p}(s,x+s-t,p).$$
We can integrate this \hbox{w.r.t.}~$p$, since $\widehat p-1$
never vanishes. Hence, the first part of the integral which appears
on the second line in~(\ref{eq:bnd:A''.1}) becomes:
\begin{eqnarray*}
I_1&:=&
\int_0^t \left(A\,\frac{\partial n}{\partial x}\right)(s,x+s-t) \,
\mathrm{d}s \\
&=&\int_0^t A(s,x+s-t) \int_{\mathbb{R}} \frac1{\widehat p-1}\,\Biggl[ 
-\frac{\partial}{\partial s}[f(s,x+s-t,p)] \\
&&\mbox{} + 
F(s,x+s-t)\,\frac{\partial f}{\partial p}(s,x+s-t,p) \Biggr]\, 
\mathrm{d}p\, \mathrm{d}s \,;
\end{eqnarray*}
then, performing integration by parts in $s$ and~$p$:
\begin{eqnarray*}
I_1&=& -A(t,x)\,\int_{\mathbb{R}} \frac{f(t,x,p)}{\widehat p-1}\, 
\mathrm{d}p + A_0(x-t)\,
\int_{\mathbb{R}} \frac{f_0(x-t,p)}{\widehat p-1}\, \mathrm{d}p \\
&&\mbox{} +
\int_0^t \frac{\partial}{\partial s}[A(s,x+s-t)]\, 
\int_{\mathbb{R}} \frac{f(s,x+s-t,p)}{\widehat p-1}\, \mathrm{d}p\, 
\mathrm{d}s \\
&& \mbox{} + \int_0^t (A\,F)(s,x+s-t)\,
\int_{\mathbb{R}} f(s,x+s-t,p)\,\frac{\mathrm{d}}{\mathrm{d}p}
\left[\frac1{\widehat p-1}\right]\, \mathrm{d}p\, \mathrm{d}s.
\end{eqnarray*}
When $p$ tends to~$+\infty$, there holds:
$$\frac1{\widehat p-1} \sim -2\,p^2,\quad
\frac{\mathrm{d}}{\mathrm{d}p}\left[\frac1{\widehat p-1}\right] =
\frac{p+\sqrt{1+p^2}}{\sqrt{1+p^2}\,\left(p-\sqrt{1+p^2}\right)} 
\sim -4\,p,$$
so that the integrations by parts mentioned above are fully justified 
using Hypotheses~\ref{hypo:f0c1},~\ref{hypo:f0},~\ref{hypo:f0:bis} and 
moreover,
\begin{eqnarray*}
\left|I_1\right| &\le& C_0 + C_1\,\|A\|_t\,\int_{\mathbb{R}} (1+p^2)\, 
f(t,x,p)\, \mathrm{d}p \\
&&\mbox{}+
C_2\,\int_0^t \left\{ \|\partial_x A\|_s + \|\partial_t A\|_s +
\|A\|_s\,\|F\|_s \right\} \int_{\mathbb{R}} (1+p^2)\, f(s,x+s-t,p)\, 
\mathrm{d}p\, \mathrm{d}s.
\end{eqnarray*}
Using Lemma~\ref{lem:moments}, we obtain a bound on the integrals 
in~$p$, so that:
$$\left|I_1\right| \le C_0 + C_1\, \Pi_3\left(t\,\|F\|_t\right)\,
\left\{ \left\|\frac{\partial A}{\partial x}\right\|_s + 
\left\|\frac{\partial A}{\partial t}\right\|_s + \|A\|_s\,(1+\|F\|_s) 
\right\}.$$
Of course, $I_2:= 
\int_0^t \left(A\,\frac{\partial n}{\partial x}\right)(s,x+s-t) 
\,\mathrm{d}s$ is bounded in the same manner. Hence~(\ref{52.3:new}).
\end{proof}

In the case of two solutions corresponding to two pairs of fields
$(E_1,A_1)$ and $(E_2,A_2)$, we have similarly:
\begin{lemma}
Let $(E_1,A_1)$ and $(E_2,A_2)$ be as in Lemma~\ref{lem:52.2:new}; and
$(E'_i,A'_i)={\cal L}(E_i,A_i)$. We have:
\begin{eqnarray}
\left\|\frac{\partial^2 A'_1}{\partial x^2}-
\frac{\partial^2 A'_2}{\partial x^2}\right\|_t 
&\le& K\, \biggl\{ \left\|\frac{\partial A_1}{\partial x} - 
\frac{\partial A_2}{\partial x} \right\|_t  +
\left\|\frac{\partial A_1}{\partial t} - 
\frac{\partial A_2}{\partial t} \right\|_t \nonumber\\
&&\mbox{} + \|A_1-A_2\|_t+ \|F_1-F_2\|_t \biggr\}
\label{56.7}
\end{eqnarray}
where the constant $K$ is polynomial in 
$\left(t,\|A_i\|_t,\|\partial_x A_i\|_t,\|\partial_t A_i\|_t,\|F_i\|_t
\right)$.
\label{lem:56.6}
\end{lemma}
\begin{proof}
This time, we have to bound the integral
$$I=\int_0^t \left( \frac{\partial(n_1\,A_1)}{\partial x} - 
\frac{\partial(n_2\,A_2)}{\partial x} \right)(s,x+s-t)\, \mathrm{d}s,$$
as well as a similar one in which the current point is $(s,x-s+t)$,
and which will be handled in exactly the same manner. Writing:
\begin{eqnarray*}
\frac{\partial(n_1\,A_1)}{\partial x} - 
\frac{\partial(n_2\,A_2)}{\partial x} &=&
A_1\, \left(\frac{\partial n_1}{\partial x} - 
\frac{\partial n_2}{\partial x} \right) +
\frac{\partial n_2}{\partial x}\, \left(A_1-A_2\right)\\ 
&&\mbox{} +
n_1\, \left(\frac{\partial A_1}{\partial x} - 
\frac{\partial A_2}{\partial x} \right) +
\frac{\partial A_2}{\partial x}\, \left(n_1-n_2\right),
\end{eqnarray*}
we split $I$ into four parts $I_1,\ I_2,\ I_3,\ I_4$, which are the 
integrals of the four terms above.
By~(\ref{44.3}) and~(\ref{eq:n1-n2:2}), we obtain:
\begin{eqnarray*}
\left|I_3\right| &\le& 
\left( C_0 + C_1\, \left\|F_2\right\|_t \right)\,t\,
\left\|\partial_x A_2 - \partial_x A_1\right\|_t \,;\\
\left|I_4\right| &\le& 
\left( C_0 + C_1\,\max_i \left\|F_i\right\|_t \right)\,
\left\|\partial_x A_2\right\|_t\, t\, \left\|F_2-F_1\right\|_t.
\end{eqnarray*}
Then $I_2$ is bounded as in the previous Lemma:
$$\left|I_2\right| \le\! C\, \Pi_3\left(t\,\|F_2\|_t\right)\,
\left\{ \|\partial_x A_1-\partial_x A_2\|_t \!+\! 
\|\partial_t A_1-\partial_t A_2\|_t 
\!+\! \|A_1-A_2\|_t\,(1+\|F_2\|_t) \right\}.$$
There remains to bound $I_1$. Performing the same computations as
in Lemma~\ref{lem:52.2:new}, we obtain:
\begin{eqnarray*}
I_1&=& -A_1(t,x)\, \int_{\mathbb{R}} 
\frac{(f_1-f_2)(t,x,p)}{\widehat p-1}\, \mathrm{d}p\\
&&\mbox{} +  \int_0^t \frac{\partial}{\partial s}[A_1(s,x+s-t)]\, 
\int_{\mathbb{R}} \frac{(f_1-f_2)(s,x+s-t,p)}{\widehat p-1}\, 
\mathrm{d}p\, \mathrm{d}s \\
&& \mbox{} + \int_0^t (A_1\,F_1)(s,x+s-t)\,
\int_{\mathbb{R}} f_1(s,x+s-t,p)\,\frac{\mathrm{d}}{\mathrm{d}p}
\left[\frac1{\widehat p-1}\right]
\, \mathrm{d}p\, \mathrm{d}s\\
&& \mbox{} - \int_0^t (A_1\,F_2)(s,x+s-t)\,
\int_{\mathbb{R}} f_2(s,x+s-t,p)\,\frac{\mathrm{d}}{\mathrm{d}p}
\left[\frac1{\widehat p-1}\right]
\, \mathrm{d}p\, \mathrm{d}s.
\end{eqnarray*}
We rearrange the last two integrals by writing, as usual, $F_1\,f_1 -
F_2\,f_2 = F_1\,(f_1-f_2) + (F_1-F_2)\,f_2$. Hence:
\begin{eqnarray*}
\left|I_1\right| &\le& C \,\biggl\{ \|A_1\|_t\, \int_{\mathbb{R}}
(1+p^2)\,(f_1-f_2)(t,x,p)\,\mathrm{d}p \\
&& \mbox{} +
\int_0^t \left\{ \|\partial_x A_1\|_s + \|\partial_t A_1\|_s +\|A_1\|_s
\|F_1\|_s \right\} \int_{\mathbb{R}} (1+p^2) (f_1-f_2)(s,x+s-t,p)\,
 \mathrm{d}p\, \mathrm{d}s\\
&&\mbox{} + \int_0^t \|A_1\|_s\, \| F_1-F_2\|_s \int_{\mathbb{R}} 
(1+p^2)\, f_2(s,x+s-t,p)\, \mathrm{d}p\, \mathrm{d}s \biggr\},
\end{eqnarray*}
which is bounded with the help of Lemmas \ref{lem:difmoments}
and~\ref{lem:moments} as:
\begin{eqnarray*}
\left|I_1\right| &\le& C\, \biggl\{
\left\{ \|\partial_x A_1\|_t + \|\partial_t A_1\|_t +\|A_1\|_t\,(1+
\|F_1\|_t) \right\} \Pi_3\left(t\,\max_i \|F_i\|_t\right)\,
\int_0^t \|F_1-F_2\|_s \\
&&\mbox{} + \|A_1\|_t\, \Pi_3\left(t\,\|F_2\|_t\right)\,
\int_0^t \|F_1-F_2\|_s.\biggr\}
\end{eqnarray*}
Putting all these bounds together, we obtain~(\ref{56.7}).
\end{proof}

Such a bound as~(\ref{56.7}) cannot be derived in the NR~case, 
basically because the velocity $\widehat p=p$ is not bounded: 
singularity formation could happen at any speed. 
The proof of the above two Lemmas does not apply, because 
neither $1/(p-1)$ nor its derivative $-1/(p-1)^2$ are locally 
integrable near $p=1$, and the integrations by parts would be 
unjustified.

\smallbreak

Nevertheless, in the QR~case we immediately deduce the following 
theorem:
\begin{theorem}
Let $(E_k,A_k)_{k\in{\mathbb{N}}}$ and $(E,A)$ be as in 
Theorem~\ref{th:cv:fort}. The sequence 
$\left(\partial_x^2 A_k\right)_{k\in{\mathbb{N}}}$ 
converges uniformly in $x$ and~$t$. 
As a consequence, the triple $(f,E,A)$ defines a characteristic 
solution to~(\ref{vla.adim}--\ref{ond.adim}) on the interval 
$\left\lbrack0,T\right\rbrack$.
\label{th:solcla:2}
\end{theorem}
\begin{proof}
The convergence of $\left(\partial_x^2 A_k\right)_{k\in{\mathbb{N}}}$
immediately follows from~(\ref{56.7}) and the convergence results 
of Theorem~\ref{th:cv:fort}, and the limit is necessarily equal to
$\partial_x^2 A$. So, $\partial_x F_k$ converges uniformly towards the 
$x$-derivative of the force given by~$(E,A)$; and we get the 
conclusion.
\end{proof}


\section{Stability of equilibria}
In this section, we restrict ourselves to classical solutions to the
system~(\ref{vla.adim}--\ref{poi.adim}).
Thus, in the NR~case, we take $T>0$ small enough to have a
classical solution on~$[0,T]$. 

Furthermore, we only consider solutions of finite mass, energy and 
entropy --- the latter concepts being precised in the following 
Subsection --- in order to rule out some unphysical or pathological 
behaviours.

\subsection{Energies and entropy for the 
system~(\ref{vla.adim}--\ref{poi.adim})}
The physical energy of the system naturally splits in two parts, 
the \emph{transversal} energy made up of the terms containing
the vector potential~$A$, and the \emph{longitudinal} energy for
the other terms.

The transversal energy is written $WT[f(t),A(t),\partial_t A(t)]$, 
where:
\begin{equation}
WT[f,A,{\dot A}] := \frac12\,\int \left\{ n[f](x)\,A(x)^2 + 
\left|{\partial_x A}(x)\right|^2 + {\dot A}(x)^2  \right\}\, 
\mathrm{d}x
\label{df:WT}
\end{equation}
and, of course, $n[f]$ denotes the density of~$f$. In the periodic 
case, the integral is to be taken over one
space period; in the open-space case, over the whole space. 

Let ${\cal V}=H^1({\mathbb{R}})$ or $H^1_L({\mathbb{R}})$ and 
${\cal H}= L^2({\mathbb{R}})$ or~$L^2_L({\mathbb{R}})$. The norm 
denoted $\|\cdot\|$, without any
subscript, will be that of~${\cal H}$; and $(\cdot\mid\cdot)$ is its 
scalar product. For any function $w\in\cesp{0}{{\cal H}}$, we set
$$\left\langle w \right\rangle_t = \sup_{0\le s\le t} \|w(s)\|^2.$$
We have the following result:
\begin{proposition}
Assume that $A_0\in{\cal V}$ and ${\dot A}_0\in{\cal H}$. Then 
$A\in\cesp{0}{{\cal V}}\cap\cesp{1}{{\cal H}}$.
As a consequence, the transversal energy is finite-valued and 
differentiable on~$[0,T]$, and there holds:
\begin{equation}
\frac{\mathrm{d}}{\mathrm{d}t} WT[f(t),A(t),\partial_t A(t)] = 
\frac12\, \int \partial_t n(t,x)\, 
A(t,x)^2\, \mathrm{d}x.
\label{eq:der:WT}
\end{equation}
\label{pro:WT}
\end{proposition}
\begin{proof}
Consider the approximating sequence 
$\left(A_k\right)_{k\in{\mathbb{N}}}$ 
from Theorem~\ref{th:cv:fort}. Clearly, $A_0\in\cesp{0}{{\cal V}}
\cap\cesp{1}{{\cal H}}$. Assuming that $A_k\in\cesp{0}{{\cal V}}\cap
\cesp{1}{{\cal H}}$, we see that $A_{k+1}$ can be identified with the
variational solution to:
\begin{eqnarray*}
&&\frac{\mathrm{d}^2}{\mathrm{d}t^2} \left( A_{k+1}(t) \mid B \right) +
\left( \partial_x A_{k+1}(t) \mid \partial_x B \right) = 
- \left( n_k(t)\,A_k(t) \mid B \right),\quad \forall B\in{\cal V},\\
&&A_{k+1}(0) = A_0\in{\cal V},\quad \partial_t A_{k+1}(0) = 
{\dot A}_0\in{\cal H},
\end{eqnarray*}
where, of course, $n_k = n[f_k]$; recall that, thanks to~\eqref{44.3},
$\|n_k\|_t$ is uniformly bounded by a constant~$n_*$ on the interval 
$[0,T]$. 
Hence $A_{k+1}\in\cesp{0}{{\cal V}}\cap\cesp{1}{{\cal H}}$.
By induction, the whole sequence belongs to this space.

Moreover, we have the classical energy estimate:
\begin{equation}
\frac{\mathrm{d}}{\mathrm{d}t}\,\frac12\,\left( \left\|
\frac{\partial A_{k+1}}{\partial t}\right\|^2 + 
\left\|\frac{\partial A_{k+1}}{\partial x}\right\|^2 \right) =
- \left( n_k(t)\,A_k(t) \Bigm| \frac{\partial A_{k+1}}{\partial t} 
\right),
\label{eq:enrg:ondes}
\end{equation}
{f}rom which we deduce, by Young's inequality:
\begin{equation}
\frac12\, \left\langle\frac{\partial A_{k+1}}{\partial t}
\right\rangle_t + \left\langle\frac{\partial A_{k+1}}{\partial x}
\right\rangle_t \le C_0 + C_1\, \left\langle A_k\right\rangle_t.
\label{a1.1}
\end{equation}
On the other hand,
\begin{displaymath}
\left\|A_k(t)\right\|^2 = \left\|A_0\right\|^2 + 2\, \int_0^t
\left( A_k(t) \Bigm| \frac{\partial A_k}{\partial t} \right)\, 
\mathrm{d}s \,; 
\end{displaymath}
applying once more Young's inequality, we get:
\begin{equation}
\left\langle A_k(t)\right\rangle_t \le C_0 + C_1\,\int_0^t 
\left\langle\frac{\partial A_k}{\partial t}\right\rangle_s\,
\mathrm{d}s.
\label{a1.2}
\end{equation}
This inequality, together with~\eqref{a1.1}, shows that
$$\frac12\, \left\langle\frac{\partial A_{k+1}}{\partial t}
\right\rangle_t + \left\langle
\frac{\partial A_{k+1}}{\partial x}\right\rangle_t \le C_0 + C_1\,
\int_0^t \left\{ \frac12\,
\left\langle\frac{\partial A_k}{\partial t}\right\rangle_s + 
\left\langle \frac{\partial A_k}{\partial x}\right\rangle_s \right\}\,
\mathrm{d}s,$$
which gives the uniform boundedness of 
$\left\langle\partial_t A_k\right\rangle_t$
and $\left\langle\partial_x A_k\right\rangle_t$ by 
Lemma~\ref{lem:56.2}, and that of $\left\langle A_k\right\rangle_t$ 
by~\eqref{a1.2}. 

As a consequence, for any~$t$ the sequences 
$\left(A_k(t)\right)_{k\in{\mathbb{N}}}$
$\left(\partial_x A_k(t)\right)_{k\in{\mathbb{N}}}$, 
$\left(\partial_t A_k(t)\right)_{k\in{\mathbb{N}}}$
admit weakly convergent subsequences in~${\cal H}$, hence 
$\left(A(t),\partial_x A(t),\partial_t A(t)\right)\in{\cal H}^3$. 
Then, similar computations show that the sequence 
$\left(A_k\right)_{k\in{\mathbb{N}}}$ is indeed strongly convergent in
$\cesp{0}{{\cal V}}\cap\cesp{1}{{\cal H}}$; and, at the 
limit~\eqref{eq:enrg:ondes} gives
\[
\frac{\mathrm{d}}{\mathrm{d}t}\,\frac12\,\left( \left\|
\frac{\partial A}{\partial t}
\right\|^2 + \left\|\frac{\partial A}{\partial x}\right\|^2 \right) =
- \left( n(t)\,A(t) \Bigm| \frac{\partial A}{\partial t} \right),
\]
from which we deduce~\eqref{eq:der:WT} by integration by parts
thanks to Corollary~\ref{cor:rg.n+j}.
\end{proof}

\medbreak

We now study the longitudinal energy. To this end, we first define the 
electrostatic field and potential functionals. In the open-space case, 
we set
\[
E[f] = -\partial_x\Phi[f]:=
\partial_x\left(\phi_\mathrm{ext} - \phi[f] \right) ,
\]
where $\phi_\mathrm{ext}$ is (the opposite of) the external confining 
potential, satisfying $\phi_\mathrm{ext}''=n_\mathrm{ext}$ and 
$\phi[f]\ge0$ is the second primitive of the density $n[f]$ of~$f$, 
given by the formula:
\begin{equation}
\phi[f](x) = \frac12\, \int_{\mathbb{R}} |x-y|\, n[f](y)\,\mathrm{d}y
\label{df:selfpot}
\end{equation}
due to Lemma \ref{lem:momentsx}.
We remark that $\phi$ is a linear, self-adjoint and positive operator.
As for the potential $\phi_\mathrm{ext}$, we assume that it satisfies 
the following hypotheses.

\begin{hypothesis}
There exist $C\ge0$ and $m\in{\mathbb{N}}$ \hbox{s.t.} 
$$
|\phi_\mathrm{ext}(x)|\le C\,\left(1+|x|^m\right) \quad\mbox{and}\quad
\int_{\mathbb{R}} |p|^m\,g(p)\,\mathrm{d}p < \infty.
$$
($g$ is the majorising function in Hypothesis~\mbox{\rm\ref{hypo:f0}}.)
\label{hypo:phiext:2}
\end{hypothesis}

Let us remark that, if $0<n_{\min}\le n_\mathrm{ext} \le n_{\max}$,
the above hypotheses are satisfied with $m=2$, provided 
Hypothesis~\ref{hypo:f0:bis} holds.

In the periodic case, we do not use the above decomposition of the
potential, since the external and the self potential would not be
periodic. We set:
\[
E[f] = -\partial_x \Phi[f],\quad -\partial_x^2 \Phi[f] = 
n_\mathrm{ext} - n[f],
\]
the uniqueness of $\Phi[f]$ being ensured by imposing periodicity and 
$$\int_0^L \Phi[f](x)\, n_\mathrm{ext}(x)\, \mathrm{d}x = 0.$$
This, in turn, ensures the self-adjointness of the operator~$\Phi$.

As for the kinetic part, we denote by $\kappa(p)$ the primitive 
of~$\widehat p$, \hbox{viz.}
\[
\kappa(p) = \frac{p^2}2 \mbox{ (NR case)}, \quad
\kappa(p) = \sqrt{1+p^2} \mbox{ (QR case).}
\]
We recall that, in the QR~case, we always assume 
Hypothesis~\ref{hypo:f0:bis}. In the NR~case, the following assumption 
will be needed to ensure the differentiability of the energy.
\begin{hypothesis}
The majorising function $g$ satisfies
$$\int_{\mathbb{R}} |p|^3\,g(p)\,\mathrm{d}p < \infty.$$
\label{hypo:f0:ter}
\end{hypothesis}
We are now ready to state the
\begin{definition}
In the periodic case, the longitudinal energy functional is given by
the two equivalent formulae:
\begin{eqnarray}
WL[f] &=& \int_0^L\! \int_{\mathbb{R}} \left\{ \kappa(p) - \frac12\,
\Phi[f](x) \right\}\, f(x,p)\, \mathrm{d}p\,\mathrm{d}x ,
\label{df:WL:perio:1}\\
WL[f] &=& \int_0^L\! \int_{\mathbb{R}} \kappa(p) \, f(x,p)\, 
\mathrm{d}p\,\mathrm{d}x + \frac12\,
\int_0^L \left|\partial_x \Phi[f](x) \right|^2\, \mathrm{d}x \,;
\label{df:WL:perio:2}
\end{eqnarray}
while in the open-space case, one sets:
\begin{equation}
WL[f] = \int_{\mathbb{R}}\int_{\mathbb{R}} \left\{ \kappa(p) +
\phi_\mathrm{ext}(x) - \frac12\,\phi[f](x) \right\}\, f(x,p)\, 
\mathrm{d}p\,\mathrm{d}x.
\label{df:WL:open}
\end{equation}
\end{definition}

Let us remark that the definition~\eqref{df:WL:open} becomes formally 
equivalent to \eqref{df:WL:perio:1} or~\eqref{df:WL:perio:2}. However,
in the open-space case $|\partial_x\phi[f]|^2$ is clearly not 
integrable and the integration by parts to arrive to a formula like 
\eqref{df:WL:perio:2} are not justified. We will comment on this in 
the next subsection. 

Now, we can analyse the evolution of the longitudinal energy.

\begin{proposition}
Let $(f,E,A)$ be a classical solution 
to~(\ref{vla.adim}--\ref{poi.adim}).
Assume Hypotheses~{\rm\ref{hypo:f0c1}}, {\rm\ref{hypo:f0}} and 
{\rm\ref{hypo:f0:ter}} if the system is~NR, 
Hypotheses~{\rm\ref{hypo:f0c1}}, {\rm\ref{hypo:f0}} and 
{\rm\ref{hypo:f0:bis}} if it is~QR. In the open-space 
setting, assume moreover that Hypotheses~{\rm\ref{hypo:f0open},
\ref{hypo:f0open:bis}} and~{\rm\ref{hypo:phiext:2}} hold.
Then $t\mapsto WL[f(t)]$ is finite-valued and differentiable for all 
$t\in[0,T]$. Namely:
\begin{equation}
\frac{\mathrm{d}}{\mathrm{d}t}WL[f(t)] = 
\int j(t,x)\,\frac{\partial}{\partial x}\frac{A^2}2
(t,x)\, \mathrm{d}x.
\label{eq:der:WL}
\end{equation}
\label{pro:WL}
\end{proposition}
\begin{proof}
First, we consider the case of periodic solutions.
The finiteness of the energy defined by~\eqref{df:WL:perio:2} 
follows from Lemma~\ref{lem:moments} and the various boundedness 
results of Theorem~\ref{th:cv:fort}. By the way, the boundedness
of~$E$ yields that of~$\Phi$; and integrating by parts gives the 
equivalent expression~\eqref{df:WL:perio:1}. 
Then, as in Corollary~\ref{cor:rg.n+j}, we bound
\begin{eqnarray*}
\frac{\partial}{\partial t} [ \kappa(p)\,f(t,x,p) ] &=& 
-\widehat p\,\kappa(p)\, \frac{\partial f}{\partial x} + 
\kappa(p)\,F(t,x)\,\frac{\partial f}{\partial p}, \\
\left| \frac{\partial}{\partial t} [ \kappa(p)\,f(t,x,p) ] \right| 
&\le& \kappa(p)\, \left( F_* + \left|\widehat p\right| \right)\, 
g_{t\,\|F\|_t}(p)\, \mathrm{e}^{t\,\left( 1+\|\partial_x F\|_t\right)}.
\end{eqnarray*}
The highest power of~$p$ is $p^3$ (NR) or $p$~(QR): the above function
is integrable in~$p$, and in~$x$ given the finite length of~$(0,L)$.
Finally, Eq.~\eqref{amp.adim} gives 
$\partial_t\left(\frac12\,E^2\right) =
E\,j$; these functions are continuous on $(0,L)$, so the second term
in~\eqref{df:WL:perio:2} is also differentiable.

Let us now consider the open-space setting. The finiteness and 
differentiablity of the integral 
$\int\!\!\!\int \left\{\kappa(p) + \phi_\mathrm{ext}(x)
\right\}\, f(t,x,p)\, \mathrm{d}x\,\mathrm{d}p$ follows from
Lemma~\ref{lem:momentsx:bis}. Then, the self potential defined 
by~\eqref{df:selfpot} satisfies:
\begin{eqnarray}
0\le \phi[f(t)](x) &\le& \frac12 \left( |x|\,\int n[f](y)\,\mathrm{d}y 
+ \int |y|\,n[f](y)\,\mathrm{d}y \right) \nonumber\\
&\le& \frac{M\,|x|}2 +\frac12\, \int\!\!\!\int (|y|+|p|)\,
f(t,y,p)\, \mathrm{d}y \le \frac12\,(M\,|x|+Z),
\label{bnd:phi}
\end{eqnarray}
$Z$ being defined in Lemma~\ref{lem:momentsx}. By the same token,
\begin{equation}
0\le \int\!\!\!\int \phi[f(t)](x)\,f(t,x,p)\,\mathrm{d}x\,\mathrm{d}p 
\le \frac12\,\int (M\,|x|+Z)\, f(t,x,p)\,\mathrm{d}x\,\mathrm{d}p 
\le M\,Z.
\label{bnd:phi*n}
\end{equation}
There remains to check the differentiability of the self 
potential energy. Using the Vlasov equation~\eqref{vla.adim} and
Corollary~\ref{cor:rg.n+j}, we find:
\begin{equation}
\left|\partial_t f(t,x,p)\right| \le \left( F_* + \left|\widehat p
\right| \right)\, g_{t\,\|F\|_t}(p)\,g_{t^2\,\|F\|_t}(x-t\,p)\,
\mathrm{e}^{t\,\left( 1+\|\partial_x F\|_t\right)}.
\label{bnd:dtf}
\end{equation}
But the linearity and self-adjointness of~$\phi$ imply: 
$\partial_t (\phi[f]\,f) =2\,\partial_t f\,\phi[f]$, which allows 
to bound
\begin{eqnarray*}
\left|\frac{\partial}{\partial t}[\phi[f(t)]\,f(t)](x,p) \right| &\le&
\left(C_0 + C_1\,|x|\right)\,
\left(1 + F_* + \left|\widehat p\right| \right) \,
g_{t\,\|F\|_t}(p)\,g_{t^2\,\|F\|_t}(x-t\,p)
\end{eqnarray*}
by~\eqref{bnd:phi} and \eqref{bnd:dtf}. This proves the integrability 
of this function, and the differentiability of the self potential 
energy.

Finally, Eq.~\eqref{eq:der:WL} is obtained through tedious but
straightforward computations (cf.~Propositions 1.5 and~1.6
in~\cite{Bouchut}), all the integrations by parts being justified
by the above arguments.
\end{proof}

From~\eqref{eq:der:WT}, \eqref{eq:der:WL} and the continuity 
equation~\eqref{eq:cons:mass}, we immediately deduce:
\begin{theorem}
Under the hypotheses of Propositions \ref{pro:WT} and~\ref{pro:WL},
the energy $W(t):= WT\left[f(t),A(t),\partial_t A(t)\right] +
WL[f(t)]$ is constant.
\label{th:cons:engy}
\end{theorem}

For reference, we give the formulae for the energy in the FR case. 
There is no splitting in transversal and longitudinal parts.
The total energy is given by:
\begin{eqnarray}
W(t) &:=& \int_0^L\!\int \left[ \sqrt{1+p^2+A(t,x)^2} - 
\frac12\,\Phi[f(t)](x) \right]\, f(t,x,p)\,\mathrm{d}p\,\mathrm{d}x 
\nonumber\\ 
&&\mbox{} +
\frac12\,\int_0^L \left[ \left|\frac{\partial A}{\partial t}(t,x)
\right|^2 + \left|\frac{\partial A}{\partial x}(t,x)\right|^2 \right]
\, \mathrm{d}x
\label{df:W:FR:perio}
\end{eqnarray}
in the periodic case, and
\begin{eqnarray}
W(t) &:=& \int_{\mathbb{R}}\int \left[ \sqrt{1+p^2+A(t,x)^2} + 
\phi_\mathrm{ext}(x) - \frac12\,\phi[f(t)](x) \right]\, 
f(t,x,p)\,\mathrm{d}p\,\mathrm{d}x \nonumber\\ 
&&\mbox{} +
\frac12\,\int_{\mathbb{R}} \left[ \left|
\frac{\partial A}{\partial t}(t,x)\right|^2 + \left|
\frac{\partial A}{\partial x}(t,x)\right|^2 \right]\, \mathrm{d}x
\label{df:W:FR:open}
\end{eqnarray}
in the open-space case. The reader may check that these expressions are
formally constant, even though the mere existence of classical 
solutions is an open problem, and subtler arguments are probably 
necessary to justify the calculations.

\begin{definition}
Let $\sigma\in{\cal C}^2((0,+\infty))\cap {\cal C}^0([0,+\infty))$ 
be a strictly convex and bounded-from-below function,
which satisfies 
$$\lim_{s\to+\infty} \frac{\sigma(s)}s = +\infty.$$
Let $\gamma$ denote the generalised inverse of~$-\sigma'$ (extended 
by~0 if necessary): it is a decreasing function in its support.

The \emph{entropy} associated to a distribution function $f$
is defined as:
\begin{displaymath}
S_\sigma[f] := \int\!\!\!\int\sigma(f(x,p))\,\mathrm{d}p\,\mathrm{d}x.
\end{displaymath}
\end{definition}

The most classical example is $\sigma(s)=s\,\ln s -s$, 
\hbox{i.e.}~$\gamma(s)=\mathrm{e}^{-s}$, associated to Maxwellian 
distribution functions.

Clearly, $S_\sigma$ is a convex, bounded-from-below, weakly lower 
semicontinuous functional on its domain of definition. 
From~\cite{CaCD02}, we know that for any $f_*\in L^1({\mathbb{R}}^d)$, 
there exists a function $\sigma$ as above \hbox{s.t.} $S_\sigma[f_*]$ 
is finite. Thus, we shall choose $\sigma$ according to the following
\begin{hypothesis}
In the periodic case, $\int \sigma(g(p))\,\mathrm{d}p < \infty$. In 
the open-space case, 
$\int \sigma(g(0)\,g(p))\,\mathrm{d}p < \infty$ and 
$\int\!\!\!\int \sigma(g(x)\,g(p))\,\mathrm{d}x\,\mathrm{d}p < \infty$.
\label{hypo:sigma}
\end{hypothesis}

In this case, the following identities hold:
\begin{eqnarray*}
\int \sigma\left(g_r(p)\right)\, \mathrm{d}p &=& 2\,r\,\sigma(g(0)) + 
\int \sigma(g(p))\,\mathrm{d}p \,; \\
\int\!\!\!\int \sigma\left(g_\rho(x)\,g_r(p)\right)\, \mathrm{d}x\,
\mathrm{d}p &=&
4\,\rho\,r\,g(0) + 2\,(\rho+r)\,\int \sigma(g(0)\,g(s))\,\mathrm{d}s \\
&& \mbox{} + \int\!\!\!\int \sigma(g(x)\,g(p))\,
\mathrm{d}x\,\mathrm{d}p.
\end{eqnarray*}

Then, the finiteness and differentiability of $S_\sigma[f(t)]$ stem 
from arguments very similar to the proof of Proposition~\ref{pro:WL}, 
and one can easily deduce.

\begin{theorem}
Under Hypothesis~{\rm\ref{hypo:sigma}},
the function $S_\sigma[f(t)]$ is constant on~$[0,T]$.
\label{th:cons:entpy}
\end{theorem}

\subsection{Equilibria of~(\ref{vla.adim}--\ref{poi.adim})}
We are looking for solutions to~(\ref{vla.adim}--\ref{poi.adim}) which
do not depend on time, \hbox{i.e.}, solutions to the coupled problem 
in $(f(x,p),A(x))$:
\begin{eqnarray}
\widehat p\, \frac{\partial f}{\partial x} + 
\frac{\mathrm{d}}{\mathrm{d}x}\left[\Phi[f] - \frac{A^2}{2}\right]\,
\frac{\partial f}{\partial p} &=& 0,
\label{vla.sta}\\
- \frac{\mathrm{d}^2 A}{\mathrm{d}x^2} + n[f]\,A &=& 0
\label{ond.sta}
\end{eqnarray}
with the potential satisfying the Poisson equation 
$-\frac{\partial^2}{\partial x^2} \Phi[f]=n_\mathrm{ext}-n[f]$.

We have the following result
\begin{lemma}
Any classical solution to~\textrm{(\ref{vla.sta}--\ref{ond.sta})} 
\hbox{s.t.}~$M>0$ and $A\in{\cal V}$ is a Vlasov--Poisson 
equilibrium, \hbox{i.e.}~it has the form~$(f_\infty,0)$, where 
$f_\infty$ solves
\begin{equation}
\widehat p\, \frac{\partial f_\infty}{\partial x} + 
\frac{\mathrm{d}\Phi[f_\infty]}{\mathrm{d}x}\,
\frac{\partial f_\infty}{\partial p} = 0, 
\label{vla.poi.sta}
\end{equation}
with $-\frac{\partial^2}{\partial x^2} \Phi[f_\infty]
=n_\mathrm{ext}-n[f_\infty]$ and consequently satisfies:
\begin{equation}
f_\infty(x,p) = {\cal F}\left(\kappa(p) - \Phi[f_\infty](x)\right),
\label{sol.vla.poi.sta}
\end{equation}
for some function~${\cal F}$.
\end{lemma}
\begin{proof}
If $A\in{\cal V}$, then for any $B\in{\cal V}$,
\begin{equation}
a_f(A,B):= \int \left( A'(x)\,B'(x) + 
n[f](x)\,A(x)\,B(x)\right) \,\mathrm{d}x = 0.
\label{eq:varia:A}
\end{equation}
Now, take $B=A$ to deduce that $a_f(A,A)=0$. Since both terms in
$a_f(A,A)$ are non-negative, we deduce that
$$
\int |A'(x)|^2 \,\mathrm{d}x = \int n[f](x)\,A(x)^2 \,\mathrm{d}x = 0.
$$
{F}rom the first identity we conclude $A$ is constant, while from the 
second we conclude $A=0$ since $n[f]>0$ on some non-negligible subset 
of ${\mathbb{R}}$ or~$(0,L)$. 
We are left with~\eqref{vla.poi.sta}, whose solution is well-known 
to be of the form~\eqref{sol.vla.poi.sta}.
\end{proof}

The function ${\cal F}$ can be precised by demanding that the solution 
should minimise a ``free energy'' functional. In other words, the 
choice of the entropy function~$\sigma$ is determined by the particular
equilibrium one is interested in.
\begin{definition}
Let $f\in L^1({\mathbb{R}}^2)$, resp.~$f\in L^1_L({\mathbb{R}}^2)$, 
the \emph{free energy} of~$f$ is:
$$K_\sigma[f]:= WL[f] + S_\sigma[f].$$
Let then $A\in {\cal V}$ and~${\dot A}\in{\cal H}$; we set
$$KT_\sigma[f,A,{\dot A}]:=  K_\sigma[f] + WT[f,A,{\dot A}].$$
\end{definition}

We consider the set of suitable distribution functions with fixed mass,
\hbox{i.e.}
\begin{displaymath}
{\cal K}(L,M) := \left\{ f\in L^1_L({\mathbb{R}}^2) : f\ge 0 
\mbox{ a.e. and } \|f\|_{L^1_L({\mathbb{R}}^2)} = M\right\}.
\end{displaymath}

\begin{lemma}
In the periodic case, $KT_\sigma$ is a strictly convex and 
bounded-from-below functional on ${\cal K}(L,M)\times{\cal V}\times
{\cal H}$. It has a unique global minimum which takes the form
$\left(f_{\infty,\sigma},0,0\right)$, where 
\begin{equation}
f_{\infty,\sigma} = \gamma\left(\kappa(p) - \Phi[f_{\infty,\sigma}](x)
-\alpha\right)
\label{eq:f.eq}
\end{equation}
and is therefore a stationary solution of the Vlasov--Poisson system.
The constant $\alpha$ is uniquely determined by $M$ and $L$.
\end{lemma}
\begin{proof}
Eq.~\eqref{df:WL:perio:1} shows that $WL$ is convex in this case; it
is clearly lower semicontinuous and on~${\cal K}(L,M)$ it is bounded 
from below. Since $S_\sigma$ enjoys the same property and is strictly 
convex, we deduce that $K_\sigma$ 
has a unique global minimum $f_{\infty,\sigma}$. Writing the Lagrange
equation expressing the minimisation under the constraint $\int f=M$
(cf.~\cite{BenAbdallah-Dolbeault01,CaCD02}), 
yields the formula~\eqref{eq:f.eq},
where $\alpha$ is the Lagrange multiplier.
Then, it is clear that for any 
$(f,A,{\dot A})\ne(f_{\infty,\sigma},0,0)$, 
$$
KT_\sigma[f,A,{\dot A}] > KT_\sigma[f,0,0] > 
KT_\sigma\left[f_{\infty,\sigma},0,0\right],
$$
where the last inequality is a consequence of the results in 
\cite{BenAbdallah-Dolbeault01,CaCD02}.
\end{proof}

A similar result in the open-space case does not hold in our case, 
in contrast to the situation studied for nonlinear stability of the 
Vlasov--Poisson system in higher dimensions.
\cite{Rein94,Braasch-Rein-Vukadinovic99,BenAbdallah-Dolbeault01,CaCD02}
The main difference in 1D being that $\partial_x \phi[f]\notin 
L^2({\mathbb{R}})$ as pointed out before. Sobolev embeddings in 
$d\geq 2$ allow, under confining conditions on the external potential 
$\phi_\mathrm{ext}$, to deduce a result similar to previous Lemma.
\cite{Dolbeault99,BenAbdallah-Dolbeault01,CaCD02} 
Moreover, the functional $KT_\sigma$ ceases to be convex in our 1D 
open-space case.

\subsection{$L^p$-nonlinear stability of equilibria in the periodic 
case}

In this subsection, we just collect the known results in several 
references and applied to the particular case we deal with.
Like in~\cite{Rein94,Braasch-Rein-Vukadinovic99,CaCD02}, one can 
rewrite:
\begin{equation}
KT_\sigma[f,A,{\dot A}] - KT_\sigma\left[f_{\infty,\sigma},0,0
\right] = \Sigma_\sigma\left[f \vert f_{\infty,\sigma}\right] + 
WT[f,A,{\dot A}],
\end{equation}
where the \emph{relative entropy} of the distribution 
function~$f$ \hbox{w.r.t.}~$g$ is defined as:
\begin{equation}
\label{relent}
\qquad \Sigma_\sigma[f\vert g] := \int\!\!\!\int 
[\sigma(f)-\sigma(g) -\sigma'(g) (f-g)] \,\mathrm{d}x\,\mathrm{d}p
+ \frac{1}{2} \int \left|\partial_x \Phi [f-g] \right|^2 \, 
\mathrm{d}x.
\end{equation}

Yet, as a consequence of Theorems~\ref{th:cons:engy} 
and~\ref{th:cons:entpy}, $KT_\sigma\left[f(t),A(t),\partial_t A(t)
\right]$ is constant for any classical solution 
to~(\ref{vla.adim}--\ref{poi.adim}).
This implies (see \cite{Rein94,Braasch-Rein-Vukadinovic99,CaCD02}) 
that:
\begin{enumerate}
\item The $L^p$ norm of $f-f_{\infty,\sigma}$ is bounded for 
$1\le p\le2$, if $\inf \sigma''(s)/s^{p-2} > 0$.
\item The $L^2$ norm of $f-f_{\infty,\sigma}$ is bounded, if 
$f_{\infty,\sigma}$ is a Maxwellian, \hbox{i.e.}~$\sigma(s) = 
s\,\ln s-s$.
\item The $H^1$ norm of $\Phi[f]-\Phi[f_{\infty,\sigma}]$ is
bounded.
\item The transversal energy is bounded --- indeed, we already knew 
this, and even a little bit more, from Proposition~\ref{pro:WT}.
\end{enumerate}
Here, all the norms are taken for $x\in [0,L]$, $p\in {\mathbb{R}}$.

The first three points follow from \S\S3 and~4 of~\cite{CaCD02} and 
references therein. Indeed, the arguments in these passages are 
independent of the dimension.

In other words, the one-dimensional periodic Vlasov--Poisson 
equilibria are $L^p$-nonlinearly stable under one-dimensional 
Vlasov--Maxwell perturbations.

Let us finally mention that Landau damping was proved in \cite{CM98}
in the case of the Vlasov--Poisson system in the periodic case. As a 
consequence, it was proved in \cite{CM98} that some Vlasov--Poisson 
equilibria which are $L^1$-nonlinearly stable, are unstable for a weak 
topology. This is not known to happen in our 1D Vlasov--Maxwell system,
although numerical computations seem to indicate that nonlinear Landau 
damping should occur in this model. Let us point out that there is no 
contradiction between these two stability assertions, since weak 
topology neighbourhoods of the 
equilibria are much larger than $L^1$ neighbourhoods.

\appendix
\section{The Duhamel formulae}
The following representation formulae are capital in the various
estimations of Sections~\ref{sec:fixptweak}--\ref{sec:clasqr}. 

The unique temperate solution~$u$ to the wave equation
\begin{displaymath}
\left\lbrace
\begin{array}{c}
{\displaystyle \frac{\partial^2 u}{\partial t^2} - 
\frac{\partial^2 u}{\partial x^2}} = f \in\lesp{1}{{\cal X}},\hfill\\
u(0,x)=u_0(x)\in{\cal X},\quad \partial_t u(0,x)=v_0(x)\in{\cal X},
\hfill
\end{array}
\right.
\end{displaymath}
where ${\cal X} = L^1_\mathrm{loc}({\mathbb{R}})\cap
{\cal S}'({\mathbb{R}})$, is explicitly given by the formula:
\begin{eqnarray}
u(t,x) &=& \frac12\,\biggl\{ u_0(x+t) + u_0(x-t) + 
\int_{x-t}^{x+t} v_0(y)\,\mathrm{d}y \nonumber\\
&&\hphantom{\frac12\,\biggl\lbrace}\mbox{} + \int_{0}^t 
\int_{x+s-t}^{x+t-s} f(s,y)\, \mathrm{d}y\,\mathrm{d}s \biggr\}.
\label{40.2}
\end{eqnarray}
Hence, after some computations, the formulae for the various 
derivatives of~$u$:
\begin{eqnarray}
\frac{\partial u}{\partial t}(t,x) &=& \frac12\, 
\biggl\{ u'_0(x+t) - u'_0(x-t) + v_0(x+t) + v_0(x-t) \nonumber\\
&&\hphantom{\frac12\,\biggl\lbrace}\mbox{}
+ \int_0^t [ f(s,x+s-t) + f(s,x+t-s) ]\,\mathrm{d}s \biggr\},
\label{40.3}\\
\frac{\partial u}{\partial x}(t,x) &=& \frac12\, 
\biggl\{ u'_0(x+t) + u'_0(x-t) + v_0(x+t) - v_0(x-t) \nonumber\\
&&\hphantom{\frac12\,\biggl\lbrace}\mbox{}
+ \int_0^t [ -f(s,x+s-t) + f(s,x+t-s) ]\,\mathrm{d}s \biggr\},
\label{40.4}\\
\frac{\partial^2 u}{\partial x^2}(t,x) &=& \frac12\,
\biggl\{ u''_0(x+t) + u''_0(x-t) + v'_0(x+t) - v'_0(x-t) \nonumber\\
&&\hphantom{\frac12\,\biggl\lbrace}\mbox{}
+ \int_0^t \left[ -\frac{\partial f}{\partial x}(s,x+s-t) + 
\frac{\partial f}{\partial x}(s,x+t-s) \right]\,\mathrm{d}s \biggr\},
\label{40.5}
\end{eqnarray}
which are valid \hbox{e.g.}~if $u'_0,\ u''_0,\ v'_0\in{\cal X}$ and
$f\in\cesp{1}{{\cal X}}$.

\section{Convergence of the sequence defined by~(\ref{df:vnt})%
\label{app:b}}
Here is the technical lemma announced in the proof of 
Theorem~\ref{th:solcla:1}.
\begin{lemma}
Let $(v_k(t))_{k\in{\mathbb{N}}}$ and $\varphi_t$ be defined 
by~(\ref{df:vnt}). 
\begin{enumerate}
\item There exists $T_1 < +\infty$ such that, for $0<t<T_1$,
$\varphi_t$ admits two fixed points $v_t < v^t$. If $v_0(t)<v^t$, 
the unstable fixed point, then the sequence $v_k(t)$ converges toward 
the stable point $v_t$.
\item For $t>T_1$, there is no fixed point, and 
$(v_k(t))_{k\in{\mathbb{N}}}$ diverges to $+\infty$.
\item When $t$ tends to zero, $v^t$ goes to infinity, while 
$v_t$ remains bounded.
\end{enumerate}
Consequently, there exists $0<T^* \le T_1$ \hbox{s.t.} 
$(v_k(t))_{k\in{\mathbb{N}}}$ is convergent for $0\le t< T^*$.
\end{lemma}
\begin{proof}
As $\varphi_t$ is convex and increasing, there
are only three possibilities:
\begin{enumerate}
\item The equation $\varphi_t(x)=x$ admits two solutions $v_t<v^t$,
which are respectively stable and unstable fixed points: 
$0<\varphi'_t(v_t)<1$ and $\varphi'_t(v^t)>1$.
\item The two solutions merge in a unique fixed point~$v$, which 
satifies both $\varphi_t(v)=v$ and $\varphi'_t(v)=1$; it is 
stable on the left side, unstable on the right side.
\item There is no fixed point at all, and $\varphi_t(x)>x$ for all
$x\in{\mathbb{R}}$.
\end{enumerate}
The case~1 is achieved for $t$ small enough. Indeed:
$$\lim_{t\to0} \beta\, t\,\mathrm{e}^{t\,(1+x)} = 0,\quad 
\forall x\in{\mathbb{R}}^+.$$
Hence, for any $\mu>1$, there exists $t_\mu>0$ \hbox{s.t.} 
$\varphi_t(\mu\,\alpha)<\mu\,\alpha$ for $t<t_\mu$. On the 
other hand, $\varphi_t(\alpha)>\alpha$. Hence, $\varphi_t$ has (at
least) one fixed point in the interval $(\alpha,\mu\,\alpha)$ 
when $t<t_\mu$. Given that $\varphi_t(x)\gg x$ when $x\to+\infty$, 
there is another fixed point in $(\mu\,\alpha,+\infty)$.

\medbreak

Then, we notice that, when $x$~is fixed, $\varphi_t(x)$ is a 
decreasing function of~$t$. This has two consequences. Firstly,
if $s<t$ and $\varphi_t$ has fixed points, then $\varphi_s$ also
has fixed points. This proves the existence of $T_1$, which may be 
finite or not. Clearly, if $T_1$ is finite, it achieves the case~2.
Eliminating~$v$ between $\varphi_t(v)=v$ and $\varphi'_t(v)=1$ gives
the following equation for~$t$:
$$\beta\,t^2\,\mathrm{e}^{\alpha\,t^2 + 2\,t} = 1.$$
As the \mbox{l.h.s.} is zero for $t=0$, infinite for $t=+\infty$ and
strictly increasing in~$t$, the equation admits a unique solution 
$T_1\in(0,+\infty)$.

\medskip

The behaviour of $v_{k+1}(t)=\varphi_t(v_k(t))$ then follows from the 
elementary theory of sequences (Figure~\ref{fig:curves}), 
and the claims 1 and~2 are obtained.

\begin{figure}[h]
\centerline{\includegraphics[width=12cm]{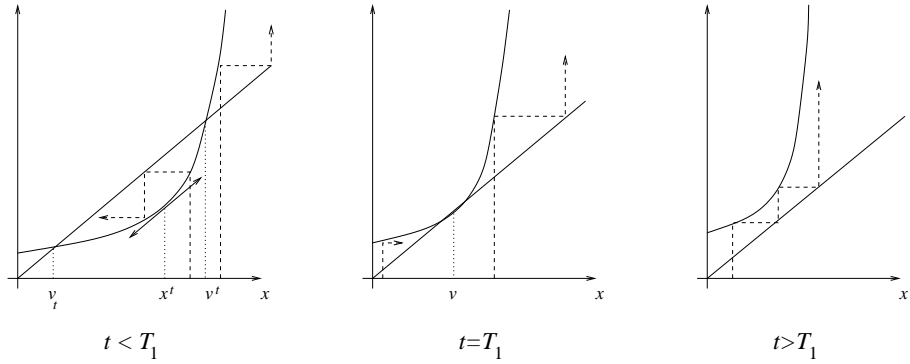}}
\caption{Dynamic of the sequence 
$v_{k+1}(t)=\varphi_t\left(v_k(t)\right)$.}
\label{fig:curves}
\end{figure}

The second consequence is that $v_t$ and $v^t$ are resp. decreasing
and increasing functions of~$t$. Hence, $v_t$ is bounded when 
$t\to0$. On the other hand, $v^t > x^t$, where $x^t$ is defined by
$\varphi'_t(x^t)=1$. The latter equation gives
$$\beta\,t^2\,\mathrm{e}^{t\,(1+x^t)}=1 \iff x^t = \frac1t\,
\ln\left(\frac1{\beta\,t^2}\right)-1 \to +\infty \quad\mbox{when } 
t\to0.$$
This proves the third claim.
Finally, let us remark that $v_0(t)$ is indeed independent of~$t$:
$E_0$ and $A_0$, hence $F_0$ are constant in time. Consequently,
for $t$ small enough, $v^t>v_0(t)$. This gives the last part of the
conclusion.
\end{proof}



\section*{Acknowledgements}
JAC and SL acknowledge support from the European IHP network
``Hyperbolic and Kinetic Equations: Asymptotics, Numerics, 
Applications'' HPRN-CT-2002-00282. JAC acknowledges the support
from the Spanish DGI-MCYT/FEDER project BFM2002-01710.
JAC thanks the hospitality of Institut \'Elie Cartan (Math\'ematiques),
Universit\'e Henri Poincar\'e Nancy during a one-month visit in which 
part of this work was achieved. We thank P.~Bertrand, T.~Goudon and 
F.~Filbet for pointing us out useful references and remarks.


\end{document}